\documentclass[12pt,leqno]{amsart}
\usepackage{amsmath,amsthm,amsfonts,amssymb}
\newtheorem{defin}{Definition}[section]
\newtheorem{lemma}[defin]{Lemma}
\newtheorem{propos}[defin]{Proposition}

\newtheorem{theor}[defin]{Theorem}
\newtheorem{corol}[defin]{Corollary}
\newtheorem{conje}[defin]{Conjecture}
\newtheorem{example}[defin]{Example}

\newenvironment{definition}{\begin{defin} \rm}{\hspace{2mm}\BOX \end{defin}}
\newenvironment{definition*}{\begin{defin} \rm}{\end{defin}}
\newcommand{\bd}{\begin{definition}}
\newcommand{\ed}{\end{definition}}
\newcommand{\bdefn}{\begin{definition*}}
\newcommand{\edefn}{\end{definition*}}
\newcommand{\bl}{\begin{lemma}}
\newcommand{\el}{\end{lemma}}
\newcommand{\bp}{\begin{propos}}
\newcommand{\ep}{\end{propos}}
\newcommand{\bt}{\begin{theor}}
\newcommand{\et}{\end{theor}}
\newcommand{\bc}{\begin{corol}}
\newcommand{\ec}{\end{corol}}
\newcommand{\bds}{\begin{displaystyle}}
\newcommand{\eds}{\end{displaystyle}}
\newcommand{\bdm}{\begin{displaymath}}
\newcommand{\edm}{\end{displaymath}}
\renewcommand{\leq}{\leqslant}
\renewcommand{\geq}{\geqslant}
\newcommand{\BOX}{\rule{2mm}{2mm}}
\newcommand{\Zint}{\mathbb Z}     % Integer number field
\newcommand{\Nat}{\mathbb N} 
\newcommand{\Rat}{\mathbb Q} 
\newcommand{\gf}{{\mathfrak g}}
\newcommand{\gl}{\mathfrak {gl}}      % Lie algebra gl
\newcommand{\sll}{\mathfrak {sl}}      % Lie algebra sl
\newcommand{\agl}{\widehat{{\mathfrak {gl}}}}      % affine Lie algebra gl
\newcommand{\asll}{\widehat{{\mathfrak {sl}}}}     % affine Lie algebra sl
\renewcommand{\ep}{\epsilon}

\newcommand{\vro}{\varrho}

\newcommand{\Vaff}{V_{\mathrm {aff}}}
\newcommand{\Vaffbar}{\ov{V}_{\mathrm {aff}}}

\newcommand{\UU}{\operatorname{U}}

\newcommand{\UNp}{\UU_q^{\prime}(\asll_N)}
\newcommand{\ULp}{\UU_q^{\prime}(\asll_L)}
\newcommand{\UM}{\UU_q(\asll_M)}
\newcommand{\UN}{\UU_q(\asll_N)}
\newcommand{\UL}{\UU_q(\asll_L)}
\newcommand{\tor}{\ddot{\UU}}
\newcommand{\Htor}{\ddot{\He}_n}
\newcommand{\EN}{E}
\newcommand{\FN}{F}
\newcommand{\KN}{K}
\newcommand{\dN}{D}
\newcommand{\EL}{\dot{E}}
\newcommand{\FL}{\dot{F}}
\newcommand{\KL}{\dot{K}}
\newcommand{\dL}{\dot{D}}
\newcommand{\Et}{\tilde{E}}
\newcommand{\Ft}{\tilde{F}}
\newcommand{\Ht}{\tilde{H}}
\newcommand{\Kt}{\tilde{K}}
\newcommand{\ct}{\tilde{c}}
\newcommand{\Xt}{\tilde{X}}
\newcommand{\Ed}{\hat{E}}
\newcommand{\Fd}{\hat{F}}
\newcommand{\Kd}{\hat{K}}
\newcommand{\pinf}{\psi_{\infty}}
\newcommand{\pinfi}{\psi_{\infty}^{-1}}
\newcommand{\pinfs}{\psi_{\infty}^2}
\newcommand{\pinfsi}{\psi_{\infty}^{-2}}
\newcommand{\mm}{\mathbf m}
\newcommand{\aab}{\mathbf a}
\newcommand{\dc}{\mathbf d}
\newcommand{\xc}{\mathbf x}
\newcommand{\M}{\operatorname M}
\newcommand{\LC}{\mathcal L}
\newcommand{\KC}{\mathcal K}
\newcommand{\EC}{\mathcal E}

\newcommand{\affHe}{\mathbf {\dot{H}}}
\newcommand{\He}{\mathbf {H}}

\newcommand{\vf}{\mathfrak v}

\newcommand{\ef}{\mathfrak e}
\newcommand{\End}{\mathrm {End}}

\newcommand{\Imm}{\mathrm {Im}}
\newcommand{\ov}{\overline}
\newcommand{\un}{\underline}
\newcommand{\unun}[1]{\un{\un{#1}}}
\newcommand{\F}{{\mathcal F}}
\newcommand{\wt}{\mathrm {wt}}

\newcommand{\Tc}{\overset{c}{T}}
\newcommand{\Ts}{\overset{s}{T}}

\newcommand{\KK}{{\mathbb K\hskip.5pt}}

\numberwithin{equation}{section}
\begin{document}

\title[Representations of the quantum toroidal algebra]{Representations of the quantum toroidal algebra on highest weight modules of the quantum affine algebra of type $\gl_N$}
\author{K. Takemura}

\address{Research Institute for Mathematical Sciences, Kyoto University, 606 Kyoto, Japan.}
\email{takemura@kurims.kyoto-u.ac.jp}
\thanks{K.T. is supported by JSPS Research Fellowship for Young Scientists} 

\author{D. Uglov}
\address{Research Institute for Mathematical Sciences, Kyoto University, 606 Kyoto, Japan.}
\email{duglov@kurims.kyoto-u.ac.jp}

%\today
\begin{abstract}
A representation of the quantum toroidal algebra of type $\sll_N$ is constructed on every integrable irreducible highest weight module of the the quantum affine algebra of type $\gl_N.$ 
The $q$-version of the level-rank duality giving the reciprocal decomposition of the $q$-Fock space with respect to mutually commutative actions of $\UU^{\prime}_q(\agl_N)$ of level $L$  and $\ULp$ of level $N$ is described.      
\end{abstract}

\maketitle

\section{Introduction}
\noindent In this article we continue our study \cite{STU} of  representations of the quantum toroidal algebra of type $\sll_N$ on irreducible integrable highest weight modules of the quantum affine algebra  of type $\gl_N.$ 
The quantum toroidal algebra $\tor$ was introduced in \cite{GKV} and \cite{VV1}. The definition of $\tor$ is given in Section \ref{sec:tor}. This algebra is a two-parameter deformation of the enveloping algebra of the universal central extension of the double-loop Lie algebra $\sll_N[x^{\pm 1}, y^{\pm 1}].$ 
To our knowledge, no general results on the representation theory of $\tor$ are available at the present. 
It therefore appears to be desirable, as a preliminary step towards a development of a general theory, to obtain concrete examples of representations of $\tor.$ 

The main reason why representations of central extensions of the double-loop Lie algebra, and of their deformations such as $\tor,$ are deemed to be a worthwhile topic to study, is that one expects  applications to higher-dimensional exactly solvable field theories. 
Our motivation to study such representations comes, however, from a different source. 
We were led to this topic while trying to understand the meaning of the level 0 action of the quantum affine algebra $\UNp$  which was defined in \cite{TU}, based on the earlier work \cite{JKKMP}, on each level 1 irreducible integrable highest weight module of the algebra $\UU_q(\agl_N).$
These level 0 actions appear as the $q$-analogues of the Yangian actions on  level 1 irreducible integrable modules of $\asll_N$ discovered in \cite{H,Schoutens}.   

Let us recall here, following \cite{STU} and \cite{VV2}, the connection between the level 0 actions and the quantum toroidal algebra $\tor.$ It is known \cite{GKV} (see also Section \ref{sec:tor}) that $\tor$ contains two subalgebras $\UU_{h},$ and $\UU_v$ such that there are algebra homomorphisms $\UNp \rightarrow \UU_h,$ and $\UNp \rightarrow \UU_v.$ 
As a consequence, every module of $\tor$ admits two actions of $\UNp:$ the {\em horizontal} action obtained through the first of the above homomorphisms, and the {\em vertical} action obtained through the second one. 
It was shown in \cite{STU} and \cite{VV2}, that on each  level 1 irreducible integrable highest weight module of $\UU_q(\agl_N)$ there is  an action of $\tor,$ such that the horizontal action coincides with the standard level 1 action of $\UNp \subset \UU_q(\agl_N),$ while the vertical action coincides with the level 0 action defined in \cite{TU}. 
The aim of the present article is to extend this result to higher level irreducible integrable highest weight modules of $\UU_q(\agl_N).$ 

The algebra $\UU_q(\agl_N)$ is, by definition, the tensor product of algebras $H\otimes\UN,$ where $H$ is the Heisenberg algebra (see Section \ref{sec:bosons}). 
Let $\Lambda$ be a level $L$ dominant integral weight of $\UN,$ and let $V(\Lambda)$ be the irreducible integrable $\UN$-module of the highest weight $\Lambda.$ As the main result of this article we define an action of $\tor$ on the irreducible $\UU_q(\agl_N)$-module    
\begin{equation}
\widetilde{V}(\Lambda) = \KK[H_-]\otimes V(\Lambda), \label{eq:tv1}
\end{equation}
where $\KK[H_-]$ is the Fock representation (see Section \ref{sec:decomp}) of $H.$ The corresponding horizontal action of $\UNp$ is just the standard, level $L,$ action on the second tensor factor in (\ref{eq:tv1}). 
The vertical action   of $\UNp$ has level zero, this action is a $q$-analogue of the Yangian action constructed recently on each irreducible integrable highest weight module of $\agl_N$ in \cite{U}.  

Let us now describe the main elements of our construction of the $\tor$-action on $\widetilde{V}(\Lambda).$ 
To define the  $\tor$-action we introduce a suitable realization of $\widetilde{V}(\Lambda)$ using the $q$-analogue of the classical level-rank duality, due to Frenkel \cite{F1,F2}, between the affine Lie algebras $\asll_N$ and $\asll_L.$ 
The quantized version of the level-rank duality takes place on the {\em $q$-Fock space} (we call it, simply, the Fock space hereafter). The Fock space is an integrable, level $L,$ module of the algebra $\UNp.$ 
The action of this algebra on the Fock space is centralized by a level $N$ action of $\ULp,$ and the resulting action of $\UNp\otimes\ULp$ is centralized by an action of the Heisenberg algebra $H.$ 

We give in the present paper a construction of the Fock space in the spirit of semi-infinite wedges of \cite{S,KMS}. The Fock space defined in \cite{KMS} appears as the special case of our construction when the level $L$ equals $1.$ 
In Theorem \ref{t:decofF} we describe the irreducible decomposition of the Fock space with respect to the action of $H\otimes\UNp\otimes\ULp.$ This theorem is the $q$-analogue of Theorem 1.6 in \cite{F1}. The decomposition shows that for every level $L$ dominant integral weight $\Lambda$ the corresponding irreducible $\UU_q(\agl_N)$-module $\widetilde{V}(\Lambda)$ is realized as a direct summand of the Fock space, such that the multiplicity space of $\widetilde{V}(\Lambda)$ is a certain level $N$ irreducible integrable highest weight module of $\ULp.$              

To define the action of the quantum toroidal algebra on $\widetilde{V}(\Lambda)$ we proceed very much along the lines of \cite{STU}. The starting point is a representation, due to Cherednik \cite{C3}, of the toroidal Hecke algebra of type $\gl_n$ on the linear space $\KK[z_1^{\pm 1},\dots,z_n^{\pm 1}]\otimes (\KK^L)^{\otimes n}.$ Here $\KK = \Rat(q^{\frac{1}{2N}}).$ 
Applying the Varagnolo--Vasserot duality \cite{VV1} between modules of the toroidal Hecke algebra and modules of $\tor,$ we obtain a representation of $\tor$ on the $q$-wedge product $\wedge^n \Vaff,$ where $\Vaff = \KK[z^{\pm1}]\otimes \KK^N\otimes \KK^L.$ This $q$-wedge product (we call it, simply, {\em the wedge product} hereafter) is similar to the wedge product of \cite{KMS}, and reduces to the latter when $L=1.$    

The Fock space is defined as an inductive limit ($n \rightarrow \infty$) of the wedge product $\wedge^n \Vaff.$ 
We show that the Fock space inherits the $\tor$-action from  $\wedge^n \Vaff.$ As the final step we demonstrate, that the $\tor$-action on the Fock space can be restricted on $\widetilde{V}(\Lambda)$ provided certain parameters in the $\tor$-action are fixed in an appropriate way.  

Let us now comment on two issues which we {\em do not} deal with in the present paper. 
The first one is the question of irreducibility of $\widetilde{V}(\Lambda)$ as the $\tor$-module. 
Based on analysis of the Yangian limit (see \cite{U}) we expect that  $\widetilde{V}(\Lambda)$ is irreducible. However we lack a complete proof of this at the present. 

The second issue is the decomposition of  $\widetilde{V}(\Lambda)$ with respect to the level 0 vertical action of $\UNp.$ 
In the Yangian limit this decomposition was performed in \cite{U} for the vacuum highest weight $\Lambda = L \Lambda_0.$ 
It is natural to expect, that combinatorially this decomposition will remain unchanged in the $q$-deformed situation. 
In particular, the irreducible components are expected to be parameterized by semi-infinite skew Young diagrams, and the $\UU_q(\sll_N)$-characters of these components are expected to be given by the corresponding skew Schur functions. \\       
\mbox{} \\
\noindent The paper is organized as follows. 
In sections \ref{s:pre} through \ref{s:Fock} we deal with the $q$-analogue of the level-rank duality, and the associated  realization of the integrable irreducible modules of $\UU_q(\agl_N).$ 
Section \ref{s:pre} contains background information on the quantum affine algebras and affine Hecke algebra. 
In Section \ref{s:wedgeprod} we introduce the wedge product, and describe the technically important {\em normal ordering rules} for the $q$-wedge vectors. In Section \ref{s:Fock} we define the Fock space, and, on this space, the action of $H\otimes\UNp\otimes\ULp.$ 
The decomposition of the Fock space as $H\otimes\UNp\otimes\ULp$-module is given in Theorem \ref{t:decofF}.    

In Sections \ref{s:tor} and \ref{s:toract} we deal with the quantum toroidal algebra $\tor$ and its actions. Section \ref{s:tor} contains basic information on the toroidal Hecke algebra and $\tor.$ 
In Section \ref{s:toract} we define actions of $\tor$ on the Fock space, and on irreducible integrable highest weight modules of $\UU_q(\agl_N).$ \\

\section{Preliminaries} \label{s:pre}
\subsection{Preliminaries on the quantum affine algebra}\label{sec:Usl}
For $k,m \in \Zint$ we define the following $q$-integers, factorials, and binomials
$$[k]_q = \frac{q^k-q^{-k}}{q-q^{-1}},\quad [k]_q!=[k]_q[k-1]_q\cdots[1]_q,\quad \text{and} \;\begin{bmatrix}m\\k\end{bmatrix}_q = \frac{[m]_q!}{[m-k]_q![k]_q!}.$$
The quantum affine algebra $\UM$ is the unital associative algebra over $\KK =\Rat(q)$ generated by the elements $E_i,F_i,K_i,K_i^{-1}, D$ $(0\leq i < M)$ subject to the relations: 
\begin{gather}
 K_i K_j = K_j K_i, \quad DK_i=K_iD,\quad K_i K_i^{-1} =  K_i^{-1} K_i = 1, \label{eq:r1}\\ 
 K_i E_j = q^{a_{ij}} E_j K_i, \label{eq:r2}\\
 K_i F_j = q^{-a_{ij}} F_j K_i, \label{eq:r3}\\ 
 [D,E_i] = \delta(i=0) E_i, \quad  [D,F_i] = -\delta(i=0) F_i,\label{eq:r4} \\
[E_i,F_j] = \delta_{ij} \frac{K_i - K_i^{-1}}{q - q^{-1}},\label{eq:r5}\\
\sum_{k=0}^{1-a_{ij}} (-1)^k \begin{bmatrix} 1- a_{ij} \\ k \end{bmatrix}_q E_i^{1-a_{ij}-k} E_j E_i^k = 0 \quad (i\neq j), \label{eq:r6}\\  
\sum_{k=0}^{1-a_{ij}} (-1)^k \begin{bmatrix} 1- a_{ij} \\ k \end{bmatrix}_q F_i^{1-a_{ij}-k} F_j F_i^k = 0 \quad (i\neq j). \label{eq:r7}
\end{gather}
Here $a_{ij} = 2\delta(i=j) - \delta(i=j+1) - \delta(i=j-1),$ and the indices are extended to all integers modulo $M.$ For $P$ a statement, we write $\delta(P) = 1$ if $P$ is true, $\delta(P) = 0$ if otherwise.    \\
$\UM$ is a Hopf algebra, in this paper we will use two different coproducts $\Delta^+$ and $\Delta^-$ given by 
\begin{alignat}{5}
&\Delta^+(K_i) = K_i\otimes K_i,\quad &  &\Delta^-(K_i) = K_i\otimes K_i,&  \label{eq:co1}\\
&\Delta^+(E_i) = E_i\otimes K_i + 1\otimes E_i,\quad & &\Delta^-(E_i) = E_i\otimes 1+ K_i\otimes E_i,& \label{eq:co2}\\   
&\Delta^+(F_i) = F_i\otimes 1 + K_i^{-1}\otimes F_i,\quad & &\Delta^-(F_i) = F_i\otimes K_i^{-1}+ 1\otimes F_i,& \label{eq:co3} \\
&\Delta^+(D) = D\otimes 1 + 1\otimes D,\quad & &\Delta^-(D) = D\otimes 1 + 1\otimes D.& \label{eq:co4}
\end{alignat}
Denote by $\UU_q^{\prime}(\asll_M)$ the subalgebra of $\UU_q(\asll_M)$ generated by $E_i,F_i,K_i,K_i^{-1}, 0\leq i < M.$ \\

\noindent In our notations concerning weights of $\UU_q(\asll_M)$ we will follow  \cite{Kac}. Thus we denote by $\Lambda_0,\Lambda_1,\dots,\Lambda_{M-1}$ the fundamental weights, by  $\delta$ the null root, and let $\alpha_i = 2\Lambda_i - \Lambda_{i+1} - \Lambda_{i-1} + \delta_{i,0} \delta$ $(0\leq i < M)$ denote the simple roots. The indices are assumed to be cyclically extended to all integers modulo $M.$ Let $P_M = \Zint \delta \oplus \left(\oplus_i \Zint \Lambda_i\right)$ be the set of integral weights. \\

\noindent Let $\KK^N$ be the $N$-dimensional vector space with basis $\vf_{1},\vf_{2},\dots,\vf_N,$ and let $\KK^L$ be the $L$-dimensional vector space with basis $\ef_{1},\ef_{2},\dots,\ef_L.$ We set $\Vaff = \KK[z^{\pm 1}]\otimes \KK^L\otimes \KK^N.$ $\Vaff$ has basis $\{ z^m \ef_{a} \vf_{\ep} \}$ where $m\in \Zint$ and $1\leq a \leq L;$ $1\leq \ep \leq N.$ Both algebras $\UN$ and $\UL$ act on $\Vaff.$ $\UN$ acts in the following way:
\begin{eqnarray}
\KN_{i}(z^m\ef_a\vf_{\ep}) &=& q^{\delta_{\ep,i} - \delta_{\ep,i+1}}z^m\ef_a\vf_{\ep}, \\  
\EN_{i}(z^m\ef_a\vf_{\ep}) &=& \delta_{\ep,i+1}z^{m+\delta_{i,0}}\ef_a\vf_{\ep-1},\\
\FN_{i}(z^m\ef_a\vf_{\ep}) &=& \delta_{\ep,i}z^{m-\delta_{i,0}} \ef_a\vf_{\ep+1},\\
\dN(z^m\ef_a\vf_{\ep}) &=& mz^m\ef_a\vf_{\ep}; 
\end{eqnarray}
where $0 \leq i < N,$ and all indices but $a$ should be read modulo $N.$ \\
The action of  $\UL$ is given by
\begin{eqnarray}
\KL_{a}(z^m\ef_b\vf_{\ep}) &=& q^{\delta_{b,L-a+1} -  \delta_{b,L-a}}z^m\ef_b\vf_{\ep}, \label{eq:ul1}\\  
\EL_{a}(z^m\ef_b\vf_{\ep}) &=& \delta_{b,L-a}z^{m+\delta_{a,0}}\ef_{b+1}\vf_{\ep},\\
\FL_{a}(z^m\ef_b\vf_{\ep}) &=& \delta_{b,L-a+1}z^{m-\delta_{a,0}}\ef_{b-1}\vf_{\ep}, \label{eq:ul2}\\
\dL(z^m\ef_a\vf_{\ep}) &=& mz^m\ef_a\vf_{\ep}.
\end{eqnarray}
where $0 \leq a< N,$ and all indices but $\ep$ are to be read modulo $L.$ Above and in what follows we put a dot over the generators of $\UL$ in order to distinguish them from the generators of $\UN.$ When both $\UN$ and $\UL$ act on the same linear space and share a vector $v$ as their weight vector, we will understand that $\wt(v)$ is a sum of weights of $\UN$ and $\UL.$ Thus 
$$ \wt(z^m\ef_a\vf_{\ep}) = \Lambda_{\ep} - \Lambda_{\ep-1} + \dot{\Lambda}_{L-a+1} - \dot{\Lambda}_{L-a} + m(\delta + \dot{\delta}).$$
Here, and from now on, we put dots over the fundamental weights, etc. of $\UL.$\\

\noindent Iterating the coproduct $\Delta^+$ (cf.(\ref{eq:co1}--\ref{eq:co4})) $n-1$ times we get an action of $\UN$ on the tensor product $\Vaff^{\otimes n}.$ Likewise 
for  $\UL,$ but in this case  we use the {\em other} coproduct $\Delta^-.$

\subsection{Preliminaries on the affine Hecke algebra} \label{sec:affH}
The affine Hecke algebra of type $\gl_n,$ $\affHe_n,$ is a unital associative algebra over $\KK$ generated by elements $T_i^{\pm 1},X_j^{\pm 1}, 1\leq i < n, 1\leq j \leq n.$ These elements satisfy the following relations:
\begin{gather}
 T_i T_i^{-1} = T_i^{-1} T_i = 1,\qquad (T_i + 1)(T_i - q^2) = 0, \label{eq:ah1}\\
T_i T_{i+1} T_i = T_{i+1} T_i  T_{i+1}, \qquad T_i T_j = T_j T_i \quad \text{if $|i-j| > 1,$} \\
 X_j X_j^{-1} = X_j^{-1} X_j = 1, \qquad X_i  X_j  = X_j  X_i, \\   
 T_i X_i T_i = q^2 X_{i+1}, \qquad T_i X_j = X_j T_i \quad \text{if $ j \neq i,i+1.$}\label{eq:ah2}
\end{gather}
The subalgebra $\He_n \subset \affHe_n$ generated by the elements $T_i^{\pm 1}$ alone is known to be isomorphic to the finite Hecke algebra of type $\gl_n.$\\

\noindent Following \cite{GRV}, \cite{KMS} we introduce a representation of $\affHe_n$ on the linear space $(\KK[z^{\pm 1}]\otimes \KK^L)^{\otimes n}.$ We will identify this space with $\KK[z_1^{\pm 1},\dots,z_n^{\pm 1}]\otimes (\KK^L)^{\otimes n}$ by the correspondence  
$$ z^{m_1}\ef_{a_1} \otimes z^{m_2}\ef_{a_2}\otimes \cdots \otimes z^{m_n}\ef_{a_n} \mapsto   z_1^{m_1}z_2^{m_2}\cdots z_n^{m_n} \otimes (\ef_{a_1} \otimes \ef_{a_2}\otimes \cdots \otimes \ef_{a_n}). $$ 
Let $E_{a,b} \in \End(\KK^L)$ be the matrix units with respect to the basis $\{\ef_a\},$ and define the {\em trigonometric R-matrix} as the following operator on $(\KK[z^{\pm 1}]\otimes \KK^L)^{\otimes 2} = \KK[z_1^{\pm 1},z_2^{\pm 1}]\otimes (\KK^L)^{\otimes 2}:$   
\begin{eqnarray*} R(z_1,z_2)& =& (q^2 z_1 - z_2)\sum_{1\leq a  \leq L} E_{a,a}\otimes E_{a,a} +  q(z_1 - z_2)\sum_{1\leq a \neq b \leq L} E_{a,a}\otimes E_{b,b} + \\ &  + &z_1(q^2 - 1)\sum_{1\leq a < b \leq L} E_{a,b}\otimes E_{b,a} + z_2(q^2 - 1)\sum_{1\leq b < a \leq L} E_{a,b}\otimes E_{b,a}.  
\end{eqnarray*}
Let $s$ be the exchange operator of factors in the tensor square $(\KK[z^{\pm 1}]\otimes \KK^L)^{\otimes 2},$ and let   
\begin{equation} \Tc _{(1,2)}:= \frac{ z_1 - q^2 z_2 }{z_1 - z_2} \cdot \left( 1 - s \cdot \frac{R(z_1,z_2)}{q^2 z_1 - z_2}\right) - 1.
\label{eq:Tc}\end{equation}
The operator $\Tc_{(1,2)}$ is known as {\em the matrix Demazure-Lusztig operator} (cf. \cite{C3}), note that it is  an element of $\End\left((\KK[z^{\pm 1}]\otimes \KK^L)^{\otimes 2}\right)$ despite the presence of the denominators in the definition. For $1 \leq i < n$ we put 
\begin{equation}
\Tc_i := 1^{\otimes (i-1)}\otimes \Tc_{(i,i+1)} \otimes 1^{\otimes (n-i-1)} \quad \in \End\left( (\KK[z^{\pm 1}]\otimes \KK^L)^{\otimes n}\right).
\label{eq:Tci}\end{equation}
\begin{propos}[\cite{C3}, \cite{GRV}, \cite{KMS}] 
The map 
\begin{equation}
X_j \mapsto z_j, \qquad T_i \mapsto \Tc_i
\end{equation}
where $z_j$ stands for the multiplication by $z_j,$ extends to a right representation of $\affHe_n$ on $(\KK[z^{\pm 1}]\otimes \KK^L)^{\otimes n}.$ 
\end{propos}
\noindent Following \cite{J} we define a left action of the finite Hecke  algebra $\He_n$ on $(\KK^N)^{\otimes n}$ by 
\begin{align}
T_i \mapsto \Ts_i :=  1^{\otimes (i-1)}\otimes \Ts \otimes 1^{\otimes (n-i-1)}, \quad \text{where $ \Ts \in \End\left((\KK^N)^{\otimes 2} \right),$} \label{eq:Tsi}\\ 
\text{and} \quad \Ts(\vf_{\ep_1}\otimes \vf_{\ep_2}) = \begin{cases} q^2 \vf_{\ep_1}\otimes \vf_{\ep_2} & \text{ if $\ep_1 = \ep_2,$} \\  
q \vf_{\ep_2}\otimes \vf_{\ep_1} & \text{ if $\ep_1 < \ep_2,$} \\
q \vf_{\ep_2}\otimes \vf_{\ep_1} + (q^2 - 1) \vf_{\ep_1}\otimes \vf_{\ep_2} & \text{ if $\ep_1 > \ep_2.$} \end{cases} \label{eq:Ts}
\end{align}

\section{The Wedge Product} \label{s:wedgeprod}

\subsection{Definition of the wedge product} \label{sec:wedge}
Identify the tensor product $\Vaff^{\otimes n}$ with \\$(\KK[z^{\pm 1}]\otimes \KK^L)^{\otimes n} \otimes (\KK^N)^{\otimes n}$ by the natural isomorphism 
$$ z^{m_1}\ef_{a_1}\vf_{\ep_1} \otimes  \cdots \otimes z^{m_n}\ef_{a_n}\vf_{\ep_n} \mapsto  \left(z^{m_1}\ef_{a_1} \otimes  \cdots \otimes z^{m_n}\ef_{a_n}\right)\otimes\left(\vf_{\ep_1} \otimes  \cdots \otimes\vf_{\ep_n}\right).   
$$
Then the operators $\Tc_i$ and $\Ts_i$ are extended on $\Vaff^{\otimes n}$ as $\Tc_i\otimes 1$ and $1 \otimes \Ts_i$ respectively. In what follows we will keep the same symbol  $\Tc_i$ to mean $\Tc_i\otimes 1,$ and likewise for $\Ts_i.$ We define the $n$-fold {\em $q$-wedge product } (or, simply, {\em the wedge product}) $\wedge^n \Vaff$  as the following quotient space:  
\begin{equation}
\wedge^n \Vaff := \Vaff^{\otimes n} / \sum_{i=1}^{n-1} \Imm (\Tc_i - \Ts_i).
\label{eq:qwp} \end{equation}
Note that under the specialization $q=1$ the operator $\Tc$ (\ref{eq:Tc}) tends to  {\em minus} the permutation operator of the tensor square $(\Rat[z^{\pm 1}]\otimes \Rat^L)^{\otimes 2},$ while the operator $\Ts$ (\ref{eq:Ts}) tends to {\em plus} the permutation operator of the tensor square $(\Rat^N)^{\otimes 2},$
so that (\ref{eq:qwp}) is a $q$-analogue of the standard exterior product.

\noindent {\bf Remark.} The wedge product is the dual, in the sense of Chari--Pressley \cite{CP} (see also \cite{C1}) of the $\affHe_n$-module $(\KK[z^{\pm 1}]\otimes \KK^L)^{\otimes n}:$ there is an evident isomorphism of linear spaces  
$$ \wedge^n\Vaff \cong (\KK[z^{\pm 1}]\otimes \KK^L)^{\otimes n}\otimes_{\He_n}(\KK^N)^{\otimes n}. $$ \\

\noindent For $ m \in \Zint_{\neq 0}$ define $B^{(n)}_m \in \End(\Vaff^{\otimes n})$ as
\begin{equation}
B^{(n)}_m = z_1^m + z_2^m + \cdots + z_n^m. 
\end{equation}
In Section \ref{sec:Usl} mutually commutative actions of the quantum affine algebras $\UN$ and $\UL$ were defined on $\Vaff^{\otimes n}.$ The operators $B^{(n)}_m$ obviously commute with these actions.\\

\noindent The following proposition is easily deduced from the results of \cite{CP}, \cite{GRV}, \cite{KMS}. 
\begin{propos} \label{p:CPD}
For each $i=1,\dots,n-1$ the subspace $\Imm(\Tc_i - \Ts_i)\subset \Vaff^{\otimes n}$ is invariant with respect to $\UN,\UL$ and $B^{(n)}_m$ $(m \in \Zint_{\neq 0}).$ Therefore  actions of $\UN,\UL$ and $B^{(n)}_m$ are defined on the wedge product $\wedge^n \Vaff.$ 
\end{propos}
\noindent It is clear that  the actions of $\UNp \subset \UN,$ $\ULp\subset\UL$ and $B^{(n)}_m$ $(m\in \Zint_{\neq 0})$ on the wedge product are mutually commutative.

\subsection{Wedges and normally ordered wedges} \label{sec:nowedges}
In the following discussion it will be convenient to relabel  elements of the basis $\{z^m\ef_a\vf_{\ep}\}$ of $\Vaff$ by single integer. We put $k = \ep-N(a+Lm)$ and denote $u_k = z^m\ef_a\vf_{\ep}.$ Then the set $\{u_k \:|\: k \in \Zint\}$ is a basis of $\Vaff.$ Let 
\begin{equation} u_{k_1}\wedge u_{k_2} \wedge \cdots  \wedge u_{k_n} \label{eq:wedge}
\end{equation} 
be the image of the tensor $u_{k_1}\otimes u_{k_2} \otimes \cdots  \otimes u_{k_n}$ under the quotient map from $\Vaff^{\otimes n}$ to $\wedge^n\Vaff.$ We will call a vector of the form (\ref{eq:wedge}) {\em a wedge} and will say that a wedge is {\em normally ordered} if $k_1>k_2>\dots>k_n.$ When $q$ is specialized to $1,$ a wedge is antisymmetric with respect to a permutation of any pair of indices $k_i,k_j,$ and the normally ordered wedges form a basis of  $\wedge^n\Vaff.$ In the general situation -- when $q$ is a parameter -- the normally ordered wedges still form a basis of $\wedge^n\Vaff.$ However the antisymmetry is replaced by a more complicated {\em normal ordering rule} which allows to express any wedge as a linear combination of normally ordered wedges.  

\mbox{}

\noindent Let us start with the case of the two-fold wedge product $\wedge^2\Vaff.$ The explicit expressions for the operators $\Tc_1$ and $\Ts_1$ lead for all $k \leq l$ to the normal ordering rule of the form  
\begin{equation}
u_k\wedge u_l = c_{kl}(q) u_l\wedge u_k + (q^2-1)\sum_{i\geq 1, l-i > k+i} c_{kl}^{(i)}(q) u_{l-i}\wedge u_{k+i}, \label{eq:norule1}
\end{equation}
were $ c_{kl}(q), c_{kl}^{(i)}(q)$ are Laurent polynomials in $q.$ In particular $c_{kk}(q)=-1,$ and thus $u_k\wedge u_k =0.$  
To describe all the coefficients in (\ref{eq:norule1}), we will employ a vector notation. For all $a,a_1,a_2 =1,\dots,L;$ $\ep,\ep_1,\ep_2 $ $=$ $1,\dots,N;$ $m_1,m_2 \in \Zint$ define the following column vectors:   
\begin{eqnarray}
X_{a,a}^{\ep_1,\ep_2}(m_1,m_2) &=&\left(\begin{array}{c} u_{\ep_1-N(a+L m_1)}\wedge u_{\ep_2-N(a+L m_2)} \\ u_{\ep_2-N(a+L m_1)}\wedge u_{\ep_1-N(a+L m_2)}\end{array}\right),\\
 Y_{a_1,a_2}^{\ep,\ep}(m_1,m_2) &=& \left(\begin{array}{c} u_{\ep-N(a_1+L m_1)}\wedge u_{\ep-N(a_2+L m_2)} \\ u_{\ep-N(a_2+L m_1)}\wedge u_{\ep-N(a_1+L m_2)}\end{array}\right), \\
 Z_{a_1,a_2}^{\ep_1,\ep_2}(m_1,m_2) &=& \left(\begin{array}{c} 
u_{\ep_1-N(a_1+L m_1)}\wedge u_{\ep_2-N(a_2+L m_2)} \\ 
u_{\ep_1-N(a_2+L m_1)}\wedge u_{\ep_2-N(a_1+L m_2)} \\ 
u_{\ep_2-N(a_1+L m_1)}\wedge u_{\ep_1-N(a_2+L m_2)} \\ 
u_{\ep_2-N(a_2+L m_1)}\wedge u_{\ep_1-N(a_1+L m_2)}  \end{array}\right).
\end{eqnarray}
Moreover let 
\begin{eqnarray}
&  X_{a,a}^{\ep_1,\ep_2}(m_1,m_2)^{\prime}=X_{a,a}^{\ep_1,\ep_2}(m_1,m_2)^{\prime\prime} = X_{a,a}^{\ep_1,\ep_2}(m_1,m_2) \quad &\text{if $ m_1\neq m_2,$} \\  
&  Y_{a_1,a_2}^{\ep,\ep}(m_1,m_2)^{\prime}=Y_{a_1,a_2}^{\ep,\ep}(m_1,m_2)^{\prime\prime} = Y_{a_1,a_2}^{\ep,\ep}(m_1,m_2) \quad &\text{if $ m_1\neq m_2,$} \\  
&  Z_{a_1,a_2}^{\ep_1,\ep_2}(m_1,m_2)^{\prime}=Z_{a_1,a_2}^{\ep_1,\ep_2}(m_1,m_2)^{\prime\prime} = Z_{a_1,a_2}^{\ep_1,\ep_2}(m_1,m_2) \quad &\text{if $ m_1\neq m_2.$}  
\end{eqnarray}
And 
\begin{eqnarray}
X_{a,a}^{\ep_1,\ep_2}(m,m)^{\prime} &=&\left(\begin{array}{c} 
0 \\ 
u_{\ep_2-N(a+L m)}\wedge u_{\ep_1-N(a+L m)}\end{array}\right),\\
Y_{a_1,a_2}^{\ep,\ep}(m,m)^{\prime} &=& \left(\begin{array}{c} 
u_{\ep-N(a_1+L m)}\wedge u_{\ep-N(a_2+L m)} \\ 
0 \end{array}\right), \\
Z_{a_1,a_2}^{\ep_1,\ep_2}(m,m)^{\prime} &=& \left(\begin{array}{c} 
u_{\ep_1-N(a_1+L m)}\wedge u_{\ep_2-N(a_2+L m)} \\ 
0 \\
u_{\ep_2-N(a_1+L m)}\wedge u_{\ep_1-N(a_2+L m)} \\ 
0 
\end{array}\right).
\end{eqnarray}
\begin{eqnarray}
X_{a,a}^{\ep_1,\ep_2}(m,m)^{\prime\prime} &=&\left(\begin{array}{c} 
u_{\ep_1-N(a+L m)}\wedge u_{\ep_2-N(a+L m)} \\ 
0
\end{array}\right),\\
 Y_{a_1,a_2}^{\ep,\ep}(m,m)^{\prime\prime} &=& \left(\begin{array}{c} 
0 \\
u_{\ep-N(a_2+L m)}\wedge u_{\ep-N(a_1+L m)}\end{array}\right), \\
 Z_{a_1,a_2}^{\ep_1,\ep_2}(m,m)^{\prime\prime} &=& \left(\begin{array}{c} 
0 \\
u_{\ep_1-N(a_2+L m)}\wedge u_{\ep_2-N(a_1+L m)} \\ 
0 \\
u_{\ep_2-N(a_2+L m)}\wedge u_{\ep_1-N(a_1+L m)}  \end{array}\right).
\end{eqnarray}
For $t \in \Zint$ introduce also the matrices:
\begin{alignat}{5}
& M_{X} = \left(\begin{array}{c c} 0 & -q \\ -q & q^2 - 1  \end{array}\right),\quad &  &M_{X}(t) = (q^2-1)\left(\begin{array}{c c} q^{2t-2} & -q^{2t-1} \\ -q^{2t-1} & q^{2t} \end{array}\right), & & \\
&  M_{Y} = \left(\begin{array}{c c} q^{-2}-1 & -q^{-1} \\ -q^{-1} & 0 \end{array}\right),\quad &  & M_{Y}(t) = (q^{-2}-1)\left(\begin{array}{c c} q^{-2t} & -q^{-2t+1} \\ -q^{-2t+1} & q^{-2t+2} \end{array}\right).& & 
\end{alignat}
\begin{gather}
  M_{Z} = \left(\begin{array}{c c c c}      
0 & 0 & -(q-q^{-1}) & -1 \\ 
0 & 0 & -1 & 0 \\
-(q-q^{-1}) & -1 & (q-q^{-1})^2 & (q-q^{-1}) \\
-1 & 0 & (q-q^{-1}) & 0 
\end{array}\right), \\
   M_{Z}(t) = \qquad \frac{q^2-1}{q^2+1} \times \\ 
  \left(\begin{array}{c c c c}     
q^{2t}-q^{-2t} & q^{2t-1}+q^{-2t+1} & -(q^{2t+1}+q^{-2t-1}) & -(q^{2t}-q^{-2t})\\
q^{2t-1}+q^{-2t+1} & q^{2t-2}-q^{-2t+2} & -(q^{2t}-q^{-2t}) & -(q^{2t-1}+q^{-2t+1}) \\
-(q^{2t+1}+q^{-2t-1}) & -(q^{2t}-q^{-2t}) & q^{2t+2}-q^{-2t-2} & q^{2t+1}+q^{-2t-1} \\
-(q^{2t}-q^{-2t}) & -(q^{2t-1}+q^{-2t+1}) & q^{2t+1}+q^{-2t-1}& q^{2t}-q^{-2t} \end{array}\right).  \nonumber 
\end{gather}
Note that all entries of the matrix $M_{Z}(t)$ are Laurent polynomials in $q,$ i.e. the numerators are divisible by  ${q^2+1}.$ 

\mbox{} 

\noindent Computing  $\Imm(\Tc - \Ts)$ we get  the  following lemma:
\begin{lemma}[Normal ordering rules] \label{l:norules}\mbox{} 
In $\wedge^2 \Vaff $ there are the following relations:
\begin{eqnarray}
 & & u_{\ep-N(a+Lm_1)}\wedge u_{\ep-N(a+Lm_2)}  =  - u_{\ep-N(a+Lm_2)}\wedge u_{\ep-N(a+Lm_1)} \quad (m_1 \geq m_2),  \label{eq:n1}
\end{eqnarray}
\begin{eqnarray}
&  &X_{a,a}^{\ep_1,\ep_2}(m_1,m_2)^{\prime}  =  M_X\cdot X_{a,a}^{\ep_1,\ep_2}(m_2,m_1)^{\prime\prime} + \!\!\!\!\sum_{t=1}^{[\frac{m_1-m_2}{2}]} \!\!\!\! M_X(t)\cdot X_{a,a}^{\ep_1,\ep_2}(m_2+t,m_1-t)^{\prime\prime}\label{eq:n2}\\ & &  \quad (m_1\geq m_2; \ep_1 > \ep_2 ), \nonumber\\
& & Y_{a_1,a_2}^{\ep,\ep}(m_1,m_2)^{\prime}  =  M_Y\cdot Y_{a_1,a_2}^{\ep,\ep}(m_2,m_1)^{\prime\prime} + \!\!\!\!\sum_{t=1}^{[\frac{m_1-m_2}{2}]}\!\!\!\! M_Y(t)\cdot Y_{a_1,a_2}^{\ep,\ep}(m_2+t,m_1-t)^{\prime\prime} \label{eq:n3}\\ & &  \quad (m_1\geq m_2; a_1 > a_2 ), \nonumber \\ 
&  &Z_{a_1,a_2}^{\ep_1,\ep_2}(m_1,m_2)^{\prime}  =  M_Z\cdot Z_{a_1,a_2}^{\ep_1,\ep_2}(m_2,m_1)^{\prime\prime} + \!\!\!\!\sum_{t=1}^{[\frac{m_1-m_2}{2}]}\!\!\!\! M_Z(t)\cdot Z_{a_1,a_2}^{\ep_1,\ep_2}(m_2+t,m_1-t)^{\prime\prime} \label{eq:n4} \\ & & \quad (m_1\geq m_2; \ep_1 > \ep_2; a_1 > a_2 ). \nonumber
\end{eqnarray}
%These relations are equivalent to {\em (\ref{eq:norule1})}. \\ 

%\noindent {\em (ii)} For $k\leq l$ define $w_{kl} \in \Vaff^{\otimes 2}$ by  
%\begin{equation}
%w_{kl} = u_k\otimes u_l - c_{kl}(q) u_l\otimes u_k - (q^2-1)\sum_{i\geq 1}^{l-i > k+i} c_{kl}^{(i)}(q) u_{l-i}\otimes u_{k+i}. \label{eq:bas}
%\end{equation}
%Here $c_{kl}(q), c_{kl}^{(i)}(q)$ are the coefficients in {\em (\ref{eq:norule1})} determined by the relations {\em (\ref{eq:n1} -- \ref{eq:n4})}.

%\noindent Then $\{w_{kl}\:|\:k \leq l\}$ is a basis of $\Imm(\Tc_1-\Ts_1) \subset  \Vaff^{\otimes 2}.$  \\

%\noindent {\em (iii)} The set of normally ordered wedges $\{ u_l \wedge u_k \: | \: l > k \}$ is a basis of $\wedge^2 \Vaff .$
\end{lemma}
\noindent The relations (\ref{eq:n1} -- \ref{eq:n4}) indeed have the form (\ref{eq:norule1}), in particular,  all wedges $u_{k}\wedge u_{l}$ in the left-hand-sides satisfy $k\leq l$ and all wedges in the right-hand-sides are normally ordered. Note moreover, that every wedge  $u_{k}\wedge u_{l}$ such that $k\leq l$ appears in the left-hand-side of one of the relations. When $L=1$ the normal ordering rules are given by (\ref{eq:n1}) and (\ref{eq:n2}), these relations coincide with the normal ordering rules of \cite[ eq.(43),(45)]{KMS}.

\begin{propos} \label{p:fwbasis}\mbox{} \mbox{}\\
{\em (i)} Any wedge from $\wedge^n \Vaff$ is a linear combination of normally ordered wedges with coefficients determined by the normal ordering rules  {\em (\ref{eq:n1} -- \ref{eq:n4})} applied in each pair of adjacent factors of  $\wedge^n \Vaff.$ \\
\noindent {\em (ii)} Normally ordered wedges form a basis of $\wedge^n \Vaff.$ 
\end{propos}
\begin{proof}
(i) follows directly from the definition of $\wedge^n \Vaff.$ \\
(ii) In view of (i) it is enough to prove that normally ordered wedges are linearly independent. This is proved by specialization $q=1.$ Let $w_{1},\dots,w_{m}$ be a set of distinct normally ordered wedges in $\wedge^n \Vaff,$ and let $t_{1},\dots,t_{m} \in  \Vaff^{\otimes n}$ be the corresponding pure tensors. Assume that 
\begin{equation}
\sum c_j(q) w_{j} = 0, \label{eq:l1}
\end{equation}
where $c_1(q),\dots,c_m(q)$ are non-zero Laurent polynomials in $q.$ Then
\begin{equation}
\sum c_j(q) t_{j} \in \sum_{i=1}^{n-1} \Imm(\Tc_i-\Ts_i).
\end{equation}
Specializing $q$ to be $1$ this gives 
\begin{equation}
\sum c_j(1) t_{j} \in \sum_{i=1}^{n-1} \Imm(P_i+1) \subset \otimes_{\Rat}^n {\Vaffbar}, \label{eq:l2}
\end{equation}
where ${\Vaffbar} = \Rat[z,z^{-1}]\otimes_{\Rat}\Rat^L\otimes_{\Rat}\Rat^N,$ and $P_i$ is the permutation operator for the $i$th and $i+1$th factors in $\otimes_{\Rat}^n {\Vaffbar}.$  Since each $t_j$ is a tensor of the form $u_{k_1}\otimes u_{k_2}\otimes \cdots \otimes u_{k_n}$ where $k_1,k_2,\dots,k_n$ is a decreasing sequence, it follows from  (\ref{eq:l2}) that $c_j(1)=0$ for all $j.$ Therefore each $c_j(q)$ has the form $(q-1)c_j(q)^{(1)}$ where $c_j(q)^{(1)}$ is a Laurent polynomial in $q.$ Equation (\ref{eq:l1}) gives now  
\begin{equation}
\sum c_j(q)^{(1)} w_{j} = 0.
\end{equation}
Repeating the arguments above we conclude that all $c_j(q)$ are divisible by arbitrarily large powers of $(q-1).$ Therefore all $c_j(q)$ vanish. 
\end{proof}

\begin{lemma}  \label{l:lemma}
Let $l \leq m.$ Then the wedges $u_m\wedge u_{m-1} \wedge \cdots \wedge u_{l+1}\wedge u_{l}\wedge u_m$ and  $u_l\wedge u_m\wedge u_{m-1} \wedge \cdots \wedge \cdots u_{l+1}\wedge u_{l}$ are equal to zero.
\end{lemma}
\begin{proof}
As  particular cases of relations (\ref{eq:n1} -- \ref{eq:n4}) we have for all  $k$ and $N\geq 2$  
$$ u_k\wedge u_k = 0, \quad u_k\wedge u_{k+1} = \begin{cases} -q^{\delta(k\not\equiv 0\bmod N)} u_{k+1}\wedge u_{k} & \text{if $N\geq 2,$}\\  -q^{-1}u_{k+1}\wedge u_{k} & \text{if $N=1.$} \end{cases} $$
The lemma follows by induction from (\ref{eq:n1} -- \ref{eq:n4}). \end{proof}

\section{The Fock Space} \label{s:Fock}

\subsection{Definition of the Fock space} 
For each integer $M$ we define the Fock space $\F_M$  as the inductive limit $(n \rightarrow \infty)$ of $\wedge^n \Vaff,$ where maps $\wedge^n \Vaff \rightarrow \wedge^{n+1} \Vaff$ are given by $v \mapsto v\wedge u_{M-n}.$ For $v \in \wedge^n \Vaff$ we denote by $v\wedge u_{M-n}\wedge u_{M-n-1}\wedge \cdots$ the image of $v$ with respect to the canonical map from $\wedge^n \Vaff$ to $\F_M.$ Note that for $v_{(n)} \in \wedge^n \Vaff,$ $v_{(r)} \in \wedge^r \Vaff,$ the equality $$  v_{(n)}\wedge u_{M-n}\wedge u_{M-n-1}\wedge \cdots  = v_{(r)}\wedge u_{M-r}\wedge u_{M-r-1}\wedge \cdots $$
holds if and only if there is $s \geq n,r$ such that 
$$  v_{(n)}\wedge u_{M-n}\wedge u_{M-n-1}\wedge \cdots \wedge u_{M-s+1} = v_{(r)}\wedge u_{M-r}\wedge u_{M-r-1}\wedge \cdots \wedge u_{M-s+1} .$$ In particular, $v_{(n)}\wedge u_{M-n}\wedge u_{M-n-1}\wedge \cdots$ vanishes if and only if there is $s \geq n$ such that $ v_{(n)}\wedge u_{M-n}\wedge u_{M-n-1}\wedge \cdots \wedge u_{M-s+1} $ is  zero. \\ 

\noindent For a decreasing sequence of integers $(k_1 > k_2 > \cdots )$ such that $k_i = M-i+1$ for $i \gg 1,$ we will call the vector 
$ u_{k_1}\wedge u_{k_2} \wedge \cdots \quad \in \F_M$ 
a (semi-infinite) {\em normally ordered wedge}.
\begin{propos}\label{p:siwbasis}
The normally ordered wedges form a basis of $\F_M.$ 
\end{propos}
\begin{proof}
For each $w \in \F_M$ there are $n, v\in \wedge^n \Vaff $ such that $ w = v \wedge u_{M-n}\wedge u_{M-n-1}\wedge \cdots.$ By Proposition \ref{p:fwbasis} the finite normally ordered wedges form a basis of $\wedge^n \Vaff,$ therefore $w$ is a linear combination of vectors 
\begin{equation}
u_{k_1}\wedge u_{k_2}\wedge \cdots \wedge u_{k_n}\wedge u_{M-n}\wedge u_{M-n-1}\wedge \cdots , \quad \text{ where $k_1 > k_2 > \cdots > k_n.$}
\label{eq:swibasis1}
\end{equation}
If $k_n \leq M-n,$ then there is $r > n$ such that  $u_{k_n}\wedge u_{M-n}\wedge u_{M-n-1}\wedge \cdots \wedge u_{M-r+1}$  vanishes by Lemma \ref{l:lemma}. It follows that (\ref{eq:swibasis1}) is zero if  $k_n \leq M-n.$ Thus the normally ordered wedges span  $\F_M.$ \\

\noindent Suppose $ \sum c_{(k_1,k_2,\dots)} u_{k_1}\wedge u_{k_2}\wedge \cdots = 0,$ where wedges under the sum are normally ordered and $c_{(k_1,k_2,\dots)} \in \KK.$ Then by definition of the inductive limit there exists  $n$ such that  $ \sum c_{(k_1,k_2,\dots)} u_{k_1}\wedge u_{k_2}\wedge \cdots \wedge u_{k_n} = 0.$ Thus 
linear independence of semi-infinite normally ordered wedges follows from the linear independence of finite normally ordered wedges.  
\end{proof}

\subsection{The actions of $\UN$ and $\UL$ on the Fock spaces} \label{sec:UF}

Define the {\em vacuum vector} of $\F_M$ as  
$$ |M\rangle  = u_M \wedge u_{M-1} \wedge \cdots . $$ 
Then for 
each vector $w$ from $\F_M$ there is a sufficiently large integer $m$ such that $w$ can be represented as 
\begin{equation}
w = v\wedge |-NLm\rangle ,\quad \text{ where $v \in \wedge^{M+NLm}\Vaff.$} \label{eq:w=vbyvac}
\end{equation} 

\noindent For each $M\in \Zint$ we define on $\F_M$ operators $\EN_i,\FN_i,\KN_i^{\pm 1}, \dN$ $(0\leq i < N)$ and  $\EL_a,\FL_a,\KL_a^{\pm 1},\dL$ $(0\leq a < L)$ and then show, in Theorem \ref{t:UNUL}, that these operators satisfy the defining relations of $\UN$ and $\UL$ respectively.

As the first step we define actions of these operators on vectors of the form $|-NLm\rangle.$ Let $\ov{v} = u_{-NLm}\wedge u_{-NLm -1}\wedge \cdots \wedge u_{-NL(m+1)+1}.$  We set 
\begin{alignat}{5}
&\dN |-NLm\rangle & = & NL\frac{m(1-m)}{2} |-NLm\rangle, \label{eq:dNvac}& \\
&\KN_i |-NLm\rangle & = & q^{L\delta(i=0)}|-NLm\rangle,\label{eq:KNvac} & \\
&\EN_i |-NLm\rangle & = & 0, \label{eq:ENvac} & \\
&\FN_i |-NLm\rangle & = & \begin{cases} 0 &  \text{ if $i \neq 0,$}  \\
 \FN_0(\ov{v})\wedge  |-NL(m+1)\rangle &  \text{ if $i = 0.$} \end{cases} & \label{eq:FNvac} 
%\end{eqnarray}
\end{alignat}
And 
\begin{alignat}{5}
&\dL |-NLm\rangle & = & NL\frac{m(1-m)}{2} |-NLm\rangle, & \label{eq:dLvac} \\
&\KL_a |-NLm\rangle & = & q^{N\delta(a=0)}|-NLm\rangle,& \label{eq:KLvac} \\
&\EL_a |-NLm\rangle & = & 0, &\label{eq:ELvac} \\
&\FL_a |-NLm\rangle & = & \begin{cases} 0 &  \text{ if $a \neq 0,$} \\
 q^{-N}\FL_0(\ov{v})\wedge  |-NL(m+1)\rangle &  \text{ if $a = 0.$} \end{cases}&\label{eq:FLvac}
\end{alignat}
Then the actions on an arbitrary vector $w \in \F_M$ are defined by using the presentation (\ref{eq:w=vbyvac}) and the coproducts (\ref{eq:co1} -- \ref{eq:co4}). Thus for $v \in   \wedge^{M+NLm}\Vaff$ and $w = v\wedge |-NLm\rangle  \in \F_M$  we define  
\begin{alignat}{5}
&\dN(w) & = & \dN(v)\wedge   |-NLm\rangle  + v \wedge \dN  |-NLm\rangle , &\label{eq:dN} \\ 
&\KN_i(w) &=& \KN_i(v)\wedge \KN_i |-NLm\rangle,&  \label{eq:KN} \\ 
&\EN_i(w) &=& \EN_i(v)\wedge \KN_i |-NLm\rangle, & \label{eq:EN} \\ 
&\FN_i(w) &=& \FN_i(v)\wedge |-NLm\rangle + \KN_i^{-1}(v)\wedge \FN_i|-NLm\rangle .&\label{eq:FN}
\end{alignat}
And 
\begin{alignat}{5}
&\dL(w) & = & \dL(v)\wedge   |-NLm\rangle  + v \wedge \dL  |-NLm\rangle , &\label{eq:dL} \\ 
&\KL_a(w) &=& \KL_a(v)\wedge \KL_a |-NLm\rangle, & \label{eq:KL} \\ 
&\EL_a(w) &=& \EL_a(v)\wedge |-NLm\rangle,  \label{eq:EL} &\\ 
&\FL_a(w) &=& \FL_a(v)\wedge \KL_a^{-1}|-NLm\rangle + v\wedge \FL_a|-NLm\rangle .&\label{eq:FL}
\end{alignat}
It follows from Lemma \ref{l:lemma} that the operators $\EN_i,\FN_i,\KN_i^{\pm 1},\dN$ and $\EL_a,\FL_a,\KL_a^{\pm 1},\dL$ are well-defined, that is do not depend on a particular choice of the presentation (\ref{eq:w=vbyvac}), and for $v \in \wedge^n \Vaff,$ $u\in \F_{M-n}$ satisfy the following relations, analogous to the coproduct formulas (\ref{eq:co1} -- \ref{eq:co4}):
\begin{alignat}{5}
&\dN(v\wedge u) & = & \dN(v)\wedge u  + v \wedge \dN(u) ,& \label{eq:codN} \\ 
&\KN_i(v\wedge u) &=& \KN_i(v)\wedge \KN_i(u), & \label{eq:coKN} \\ 
&\EN_i(v\wedge u) &=& \EN_i(v)\wedge \KN_i (u) + v\wedge \EN_i(u), & \label{eq:coEN} \\ 
&\FN_i(v\wedge u) &=& \FN_i(v)\wedge u  + \KN_i^{-1}(v)\wedge \FN_i(u) .&\label{eq:coFN}
\end{alignat}
And 
\begin{alignat}{5}
&\dL(v\wedge u) & = & \dL(v)\wedge u  + v \wedge \dL(u) ,& \label{eq:codL} \\ 
&\KL_a(v\wedge u) &=& \KL_a(v)\wedge \KL_a(u), & \label{eq:coKL} \\ 
&\EL_a(v\wedge u) &=& \EL_a(v)\wedge u + \KL_a(v)\wedge \EL_a(u), & \label{eq:coEL} \\ 
&\FL_a(v\wedge u) &=& \FL_a(v)\wedge  \KL_a^{-1}(u)  + v\wedge \FL_a(u) .\label{eq:coFL}&
\end{alignat}
Relations ((\ref{eq:dNvac}, \ref{eq:KNvac}),(\ref{eq:dN}, \ref{eq:KN})) and  ((\ref{eq:dLvac}, \ref{eq:KLvac}),(\ref{eq:dL}, \ref{eq:KL})) define the weight decomposition of the Fock space $\F_M.$ We have   
\begin{equation} \wt( |-NLm\rangle ) = L\Lambda_0 + N\dot{\Lambda}_0 + NL\frac{m(1-m)}{2} (\delta + \dot{\delta}), \label{eq:wt1} \end{equation} 
and for $v \in \wedge^{M+NLm} \Vaff$
\begin{equation} \wt(v\wedge |-NLm\rangle ) = \wt(v) + \wt( |-NLm\rangle ). \label{eq:wt2} \end{equation} 

\begin{theor} \label{t:UNUL} \mbox{} \\
{\em (i)} The operators $\EN_i,\FN_i,\KN_i,\dN$ $(0\leq i <N)$ define on $\F_M$ a structure of an integrable $\UN$-module. And the operators  $\EL_a,\FL_a,\KL_a,\dL $ define on $\F_M$ a structure of an integrable $\UL$-module. \\
{\em (ii)} The actions of the subalgebras $\UNp \subset \UN$ and $\ULp \subset \UL$ on $\F_M$ are mutually commutative.
\end{theor}
\begin{proof}
(i) It is straightforward to verify that the relations (\ref{eq:r1}--\ref{eq:r4}) are satisfied. In particular, the weights of $\EN_i,\FN_i$ and $\EL_a,\FL_a$ are $\alpha_i,-\alpha_i$ and  $\dot{\alpha}_a,-\dot{\alpha}_a$ respectively.
To prove the relations  
\begin{equation}
[\EN_i,\FN_j] = \delta_{ij}\frac{\KN_i-\KN_i^{-1}}{q -q^{-1}}, \quad \text{and} \quad [\EL_a,\FL_b] = \delta_{ab}\frac{\KL_a-\KL_a^{-1}}{q -q^{-1}} \label{eq:t1}
\end{equation}
it is enough, by (\ref{eq:KN}--\ref{eq:FN}) and (\ref{eq:KL}--\ref{eq:FL}), to show that these relations hold when applied to a vacuum vector of the form $|-NLm\rangle.$ If $i\neq j ,$ $a\neq b$ we  have   
$$ [\EN_i,\FN_j]|-NLm\rangle = 0,\quad  [\EL_a,\FL_b]|-NLm\rangle = 0 $$
because $\alpha_i-\alpha_j + \wt(|-NLm\rangle)$ $(i\neq j )$ and $\dot{\alpha}_a -\dot{\alpha}_b+\wt(|-NLm\rangle)$ $(a\neq b)$ are not weights of $\F_{-NLm}.$
The relations 
\begin{alignat}{4} &[\EN_i,\FN_i]|-NLm\rangle &=&\frac{\KN_i-\KN_i^{-1}}{q -q^{-1}}|-NLm\rangle   \label{eq:t21} \\ &[\EL_a,\FL_a] |-NLm\rangle &= &\frac{\KL_a-\KN_a^{-1}}{q -q^{-1}}|-NLm\rangle  & \label{eq:t22}
\end{alignat}evidently hold by (\ref{eq:KNvac} -- \ref{eq:FNvac}), (\ref{eq:KLvac} -- \ref{eq:FLvac}) when $i \neq 0,$ $a\neq 0.$ 
Let $a=0.$ We have 
\begin{align*} \FL_0 |-NLm\rangle  = &q^{-N}\sum_{i =1}^N q^{i} \: u_{N-N(1+Lm)}\wedge u_{N-1-N(1+Lm)}\wedge \cdots \\ & \cdots \wedge u_{i-N(L+L(m-1))} \wedge \cdots  \wedge u_{1-N(1+Lm)} \wedge |N-N(2+Lm)\rangle.
\end{align*}
Then by Lemma \ref{l:lemma} 
$$ \EL_0 \FL_0 |-NLm\rangle = q^{1-N} \sum_{i=1}^N q^{2(i-1)} |-NLm\rangle = \frac{q^N-q^{-N}}{q -q^{-1}}|-NLm\rangle.$$  
This shows the  relation  (\ref{eq:t22}) for $a=0.$ The  relation (\ref{eq:t21}) for  $i=0$ is shown in a similar way.

 Thus $\EN_i,\FN_i,\KN_i,\dN$ and $\EL_a,\FL_a,\KL_a,\dL $  satisfy the defining relations (\ref{eq:r1} -- \ref{eq:r5}). Observe that for $i=0,\dots,N-1;$ $a=0,\dots,L-1$ and $\mu \in P_N+P_L,$ $\mu  + r \alpha_i$, $ \mu + n \dot{\alpha}_a $ are weights of $\F_{M}$ for only a finite number of $r$ and $n.$ Therefore $\F_M$ is an integrable module of $\UU_q(\sll_2)_i = \langle \EN_i,\FN_i,\KN_i^{\pm 1}\rangle $ and $\UU_q(\sll_2)_a = \langle \EL_a,\FL_a,\KL_a^{\pm 1}\rangle .$ By Proposition B.1 of \cite{KMPY} this implies that the Serre relations (\ref{eq:r6}, \ref{eq:r7}) are satisfied. 

Eigenspaces of the operator $\dN$ and eigenspaces of the operator $\dL$  are finite-dimensional. Therefore the integrability with respect to each  $\UU_q(\sll_2)_i$ and $\UU_q(\sll_2)_a$ implies the integrability of $\F_M$ as both $\UN$-module and $\UL$-module.  \\  

\noindent (ii) The Cartan part of $\UNp$ evidently commutes with $\ULp,$ and vice-versa. By (\ref{eq:KN} -- \ref{eq:FN}) and (\ref{eq:KL} -- \ref{eq:FL}) it is enough to prove that commutators between the other generators vanish when applied to a vector of the form $|-NLm\rangle .$ 
The relation 
$$ [\EN_i,\EL_a] |-NLm\rangle  =  0$$
is trivially satisfied by (\ref{eq:ENvac}, \ref{eq:ELvac}). The relations 
$$ [\FN_i,\EL_a] |-NLm\rangle  =  0, \quad [\EN_i,\FL_a] |-NLm\rangle  =  0$$
hold because $\dot{\alpha}_a - \alpha_i + \wt(|-NLm\rangle)$ and $ \alpha_i -\dot{\alpha}_a  + \wt(|-NLm\rangle)$ are not weights of $\F_{-NLm}.$
The relations $$[\FN_i,\FL_a] |-NLm\rangle  =  0$$ are trivial  by (\ref{eq:FNvac}, \ref{eq:FLvac}) when $i\neq 0,a\neq 0;$ and are verified  by using the normal ordering rules (\ref{eq:n1} -- \ref{eq:n4}) and  Lemma \ref{l:lemma} in the rest of the cases. 
\end{proof}

\subsection{The actions of Bosons} \label{sec:bosons}
We will now define  actions of  operators $B_n$ $(n\in \Zint_{\neq 0})$ (called {\em bosons}) on $\F_M.$  Let $u_{k_1}\wedge u_{k_2} \wedge \cdots $ ($k_i = M-i+1$ for $i\gg 1$) be a vector of $\F_M.$ By Lemma \ref{l:lemma}, for $n\neq 0$ the sum 
\begin{align}
 &(z^n u_{k_1})\wedge u_{k_2} \wedge u_{k_3} \wedge  \cdots \; + \label{eq:ba}\\ 
 &u_{k_1} \wedge(z^n u_{k_2})\wedge u_{k_3} \wedge  \cdots \; + \nonumber\\ 
 &u_{k_1} \wedge u_{k_2} \wedge(z^n u_{k_3})\wedge \cdots \; + \nonumber\\ 
  & \quad + \quad \cdots \quad .\nonumber 
\end{align}
contains only a finite number of non-zero terms, and is, therefore, a vector of $\F_M.$ By Proposition \ref{p:CPD} the assignment $u_{k_1}\wedge u_{k_2} \wedge \cdots  \mapsto \text{(\ref{eq:ba})}$ defines an operator on $\F_M.$ We denote this operator $B_n.$ By definition we have for $v \in \Vaff,$ $u \in \F_{M-1}:$
\begin{equation}
 B_n(v\wedge u) = (z^nv)\wedge u + v\wedge B_n(u). \label{eq:Bvu}
\end{equation}
\begin{propos}
For all $n\in \Zint_{\neq 0}$ the operator $B_n$ commutes with the actions of $\UNp$ and $\ULp.$ 
\end{propos}
\begin{proof}
It follows immediately from the definition, that the weight of $B_n$ is $n(\delta + \dot{\delta}).$ Thus  $B_n$ commutes with $\KN_i, \KL_a$ $(0\leq i <N, 0\leq a < L).$ 

Let $X$ be any of the operators $\EN_i,\FN_i,\EL_a,\FL_a$ $(0\leq i <N, 0\leq a < L).$  The relations (\ref{eq:coEN}, \ref{eq:coFN}), (\ref{eq:coEL}, \ref{eq:coFL}) and (\ref{eq:Bvu})  imply now that $[B_n , X] = 0$ will follow from  $[B_n , X] |-NLm \rangle = 0 $ for an arbitrary integer $m.$    

If $n >0,$ we have $[B_n , X] |-NLm \rangle = 0$ because $n(\delta + \dot{\delta})\pm\alpha_i + \wt(|-NLm \rangle) $ and $n(\delta + \dot{\delta})\pm\dot{\alpha}_a + \wt(|-NLm \rangle)$ are not weights of $\F_{-NLm}.$  

Let $n < 0.$ Consider the expansion 
$$
[B_n , X] |-NLm \rangle  = \sum_{\nu} c_{\nu} u_{k^{\nu}_1}\wedge u_{k^{\nu}_2}\wedge \cdots 
$$
where the wedges in the right-hand-side are normally ordered. Comparing the weights of the both sides, we obtain  for all $\nu$ the inequality $k_1^{\nu} > -NLm.$ For $r \geq 0$  (\ref{eq:coEN}, \ref{eq:coFN}), (\ref{eq:coEL}, \ref{eq:coFL}) and (\ref{eq:Bvu}) give 
\begin{multline}
[B_n , X] |-NLm \rangle  = \\ = u_{-NLm}\wedge u_{-NLm-1}\wedge \cdots \wedge u_{-NL(m+r)+1}\wedge [B_n , X] |-NL(m+r) \rangle   \label{eq:bx}
\end{multline}
where 
$$ [B_n , X] |-NL(m+r) \rangle   = 
\sum_{\nu} c_{\nu} u_{k^{\nu}_1-NLr}\wedge u_{k^{\nu}_2-NLr}\wedge \cdots .
$$
Now let $r$ be sufficiently large, so that 
$$ k_1^{\nu} - NLr \leq   -NLm$$
holds for all $\nu.$ By Lemma \ref{l:lemma}, the last inequality and $ k_1^{\nu} - NLr > -NL(m+r) $ imply that (\ref{eq:bx}) vanishes. 
\end{proof}

\begin{propos}
There are non-zero $\gamma_n(q) \in \Rat[q,q^{-1}]$ (independent on $M$) such that 
\begin{equation}
[B_n,B_{n^{\prime}}] = \delta_{n+n^{\prime},0} \gamma_n(q).
\end{equation}
\end{propos}
\begin{proof}
Each vector of $\F_{M^{\prime}}$ $(M^{\prime}\in\Zint)$ is of the form $v\wedge |M\rangle$ where $v \in \wedge^k \Vaff,$ and $k=M^{\prime}-M$ is sufficiently large. By (\ref{eq:Bvu}) we have  
$$ [B_n,B_{n^{\prime}}](v\wedge |M\rangle) = v\wedge[B_n,B_{n^{\prime}}]|M\rangle.$$
The vector $[B_n,B_{n^{\prime}}]|M\rangle$ vanishes if $n+n^{\prime}>0$ because in this case  $\wt(|M\rangle) + (n+n^{\prime})(\delta + \dot{\delta})$ is not a weight of $\F_M.$

Let $n+n^{\prime}<0.$ Write $[B_n,B_{n^{\prime}}]|M\rangle$ as the linear combination of normally ordered wedges:
$$ [B_n,B_{n^{\prime}}]|M\rangle = \sum_{\nu} c_{\nu} u_{k_1^{\nu}}\wedge u_{k_2^{\nu}} \wedge \cdots .$$ 
Since $[B_n,B_{n^{\prime}}]|M\rangle$ is of the weight $\wt(|M\rangle) + (n+n^{\prime})(\delta + \dot{\delta})$ with $n+n^{\prime}<0,$ we necessarily have $k_1^{\nu} > M.$ 
For any $s > 0$ eq. (\ref{eq:Bvu}) gives    
\begin{equation} [B_n,B_{n^{\prime}}]|M\rangle = u_M\wedge u_{M-1} \wedge \cdots \wedge u_{M-NLs+1}\wedge [B_n,B_{n^{\prime}}]|M-NLs\rangle, 
\label{eq:pB1}\end{equation}
where 
$$ 
 [B_n,B_{n^{\prime}}]|M-NLs\rangle = \sum_{\nu} c_{\nu} u_{k_1^{\nu}-NLs}\wedge u_{k_2^{\nu}-NLs} \wedge \cdots.
$$
Taking $s$ sufficiently large so that $ M-k_1^{\nu} + NLs \geq 0$ holds for all $\nu$ above, we have for all $\nu$ the inequalities 
$$ k_1^{\nu} - NLs - (M - NLs) > 0, \quad \text{and} \quad M-(k_1^{\nu} - NLs) \geq 0. $$ 
Lemma \ref{l:lemma} now shows that (\ref{eq:pB1}) is zero.

Let now $n+n^{\prime}=0.$ The vector  $[B_n,B_{n^{\prime}}]|M\rangle$ has weight $\wt(|M\rangle).$ The weight subspace of this weight is one-dimensional, so we have $ [B_n,B_{-n}]|M\rangle  = \gamma_{n,M}(q) |M\rangle $ for $\gamma_{n,M}(q) \in \KK.$ Since $[B_n,B_{-n}]|M\rangle = u_M \wedge [B_n,B_{-n}]|M-1\rangle,$  $ \gamma_{n,M}(q)$ is independent on $M.$  

The coefficients $c_{kl}(q),c_{kl}^{(i)}(q)$ in the normal ordering rules (\ref{eq:norule1}) are Laurent polynomials in $q,$ hence so are $\gamma_{n}(q).$ Specializing to $q=1$ we have $\gamma_{n}(1) = nNL.$ Thus all  $\gamma_{n}(q)$ $(n \in \Zint_{\neq 0})$ are non-zero.
\end{proof}

\begin{propos}
If $N=1$ or $L=1$ or $n=1,2$, we have for $\gamma _n(q)$ the following formula:
\begin{equation}
\gamma _n(q) = n \frac{1-q^{2Nn}}{1-q^{2n}} \frac{1-q^{-2Ln}}{1-q^{-2n}}.
\label{Bconst}
\end{equation}
\end{propos}
\begin{proof}
The $L=1$ case is due to \cite{KMS}, and the formula for $N=1$ is obtained from  the formula for $L=1$ by  comparing the normal ordering rules (\ref{eq:n2}) and (\ref{eq:n3}).
The $n=1,2$ case is shown by a direct but lengthy calculation. 
(First act with $B_{-n}$ on the vacuum vector, express all terms as linear combinations of the  normally  ordered wedges, then act with $B_n$ and, again,  rewrite the result in terms of the  normally  ordered wedges to  get the coefficient $\gamma _n(q)$.)
\end{proof}

\begin{conje}
The formula {\em (\ref{Bconst})} is valid for all positive integers $N,L,n$.
\end{conje}

\noindent Let $H$ be the Heisenberg algebra generated by $\{ B_n\}_{n\in\Zint_{\neq 0}}$ with the defining relations $[B_n,B_{n^{\prime}}]=\delta_{n+n^{\prime},0}\gamma_n(q).$ Summarizing this and the previous sections, we have constructed on each Fock space $\F_M$ an action of the algebra $H\otimes \UNp \otimes \ULp.$ Note that the action of $\UNp$ has level $L$ and the action of $\ULp$ has level $N.$	

\subsection{The decomposition of the Fock space} \label{sec:decomp}
Let $P_N^+$ and $P_N^+(L)$ be respectively the set of dominant integral weights of $\UNp$ and the subset of dominant integral weights of level $L\in \Nat:$  
\begin{alignat}{4}
&P_N^+ &= &\{ a_0\Lambda_0 + a_1\Lambda_1 + \cdots  + a_{N-1}\Lambda_{N-1} \: | \: a_i \in \Zint_{\geq 0} \}, \\
&P_N^+(L) &= &\{ a_0\Lambda_0 + a_1\Lambda_1 + \cdots  + a_{N-1}\Lambda_{N-1} \: | \: a_i \in \Zint_{\geq 0},\; \sum a_i = L \}. 
\end{alignat}
For $\Lambda \in P_N^+$ let $V(\Lambda)$ be the irreducible integrable highest weight module of $\UNp,$ and let $v_{\Lambda} \in V(\Lambda)$ be the highest weight vector.  

Let $\ov{\Lambda}_1,\ov{\Lambda}_2,\dots,\ov{\Lambda}_{N-1}$ be the fundamental weights of $\sll_N,$ and let $\ov{\alpha}_i = 2 \ov{\Lambda}_i-\ov{\Lambda}_{i+1}-\ov{\Lambda}_{i-1}$ $ 1 \leq i < N$ be the simple roots. Here the indices are cyclically extended to all integers modulo $N,$ and $\ov{\Lambda}_0 :=0.$ Let $\ov{Q}_N = \oplus_{i=1}^{N-1} \Zint \ov{\alpha}_i $ be the root lattice of $\sll_N.$ For an $\UNp$-weight $\Lambda = \sum_{i=0}^{N-1} a_i \Lambda_i $ we will set  $\ov{\Lambda} = \sum_{i=1}^{N-1} a_i \ov{\Lambda}_i.$ 

A vector $w \in \F_M$ is a {\em highest weight vector} of $H\otimes \UNp \otimes \ULp$ if it is a highest weight  vector with respect to $\UNp$ and $\ULp$ and is annihilated by $B_n$ with $n > 0.$ We will now describe a family of  highest weight vectors.

With every  $\Lambda = \sum_{i=0}^{N-1}a_i \Lambda_i \in P_N^+(L),$ such that $\ov{\Lambda} \equiv \ov{\Lambda}_M \bmod \ov{Q}_N,$ we associate $\dot{\Lambda}^{(M)} \in P_L^+(N)$ (i.e. $\dot{\Lambda}^{(M)}$ is a dominant integral weight of $\ULp$ of level $N$) as follows. Let $M\equiv s \bmod NL$ $( 0\leq s < NL),$ and let $l_1 \geq l_2 \geq \dots \geq l_N $ be the partition defined by the relations:
\begin{align}
&l_i - l_{i+1} = a_i\quad (1 \leq i < N), \label{eq:part1}\\ 
& l_1 + l_2 + \cdots + l_N = s + NL.       \label{eq:part2}
\end{align}
Note that all $l_i$ are integers, and that $l_N > 0.$ Then we set  
\begin{equation}
\dot{\Lambda}^{(M)} := \dot{\Lambda}_{l_1} + \dot{\Lambda}_{l_2} + \cdots + \dot{\Lambda}_{l_N}. \label{eq:dotLambda}  
\end{equation}
Recall that the indices of the fundamental weights are cyclically extended to all integers modulo $L$. 
Consider the Young diagram of $l_1 \geq l_2 \geq \dots \geq l_N$ (Fig. 1).We set the coordinates $(x,y)$ of the lowest leftmost square to be $(1,1).$ 
\begin{center}
\begin{picture}(120,120)(-15,-15)
\multiput(0,0)(10,0){2}{\line(0,1){85}} 
\put(0,85){\line(1,0){10}} \put(0,85){\makebox(12,15){{\scriptsize $l_1$}}} 
\put(20,0){\line(0,1){60}} 
\put(10,60){\line(1,0){10}} \put(10,60){\makebox(12,15){{\scriptsize $l_2$}}}
\put(30,0){\line(0,1){40}} 
\put(20,40){\line(1,0){10}} \put(20,40){\makebox(12,15){{\scriptsize $l_3$}}}
\put(40,15){\makebox(30,10){$\cdots$}}
\multiput(80,0)(10,0){2}{\line(0,1){20}} 
\put(80,20){\line(1,0){10}} \put(80,20){\makebox(12,15){{\scriptsize $l_N$}}}
\put(0,0){\line(1,0){30}} \put(80,0){\line(1,0){10}}

%\put(0,74){\makebox(11,15){{\scriptsize $1$}}}
%\put(0,66){\makebox(11,15){{\scriptsize $2$}}}
%\put(0,58){\makebox(11,15){{\scriptsize $3$}}}
%\put(0,50){\makebox(11,15){{\scriptsize $5$}}}
%\put(10,50){\makebox(11,15){{\scriptsize $4$}}}
%\put(0,42){\makebox(11,15){{\scriptsize $7$}}}
%\put(10,42){\makebox(11,15){{\scriptsize $6$}}}

\put(-10,-10){\vector(1,0){20}} \put(10,-15){\makebox{{\scriptsize $x$}}}
\put(-10,-10){\vector(0,1){20}} \put(-15,10){\makebox{{\scriptsize $y$}}}
\end{picture}

{\large { Fig. 1}}  \end{center}
Introduce a  numbering of squares of the Young diagram by $1,2,\dots,s+NL$ by requiring that the numbers assigned to squares in the bottom row of a pair of any adjacent rows are greater than the numbers assigned to squares in the top row, and that the numbers increase from right to left within each row (cf. the example below). Letting $(x_i,y_i)$ to be the coordinates of the $i$th square, set $k_i = x_i + N(y_i - L - 1) + M-s.$ Then $k_i > k_{i+1}$ for all $i=1,2,\dots,s+NL-1.$ Now define  
\begin{equation}
\psi_{\Lambda} = u_{k_1}\wedge u_{k_2} \wedge \cdots \wedge u_{k_{s+NL}}\wedge |M-s-NL \rangle. \label{eq:hwv}
\end{equation}
Note that $\psi_{\Lambda} \in \F_M,$ and $\psi_{\Lambda}$ is a normally ordered wedge.
\begin{example}
Let $N=3,$ $L=2,$ and $M=0.$ The set $\{ \Lambda \in P_3^+(2) \: | \: \ov{\Lambda} \equiv 0\bmod \ov{Q}_3 \}$ contains the two weights: $2\Lambda_0$ and $\Lambda_1+\Lambda_2$ only. The corresponding weights of $\ULp$ and the numbered Young diagrams  are shown below.  
$$ \begin{picture}(100,60)(0,0) 
\put(0,50){$\Lambda = 2\Lambda_0:$}
\put(0,35){$\dot{\Lambda}^{(0)} = 3\dot{\Lambda}_0$}
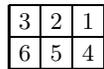
\multiput(0,0)(0,10){3}{\line(1,0){30}}
\multiput(0,0)(10,0){4}{\line(0,1){20}}
\put(0,0){\makebox(10,10){\scriptsize{$6$}}}
\put(10,0){\makebox(10,10){\scriptsize{$5$}}}
\put(20,0){\makebox(10,10){\scriptsize{$4$}}}
\put(0,10){\makebox(10,10){\scriptsize{$3$}}}
\put(10,10){\makebox(10,10){\scriptsize{$2$}}}
\put(20,10){\makebox(10,10){\scriptsize{$1$}}}
\end{picture}
\begin{picture}(100,60)(0,0) 
\put(0,50){$\Lambda = \Lambda_1+\Lambda_2:$}
\put(0,35){$\dot{\Lambda}^{(0)} = \dot{\Lambda}_0 + 2\dot{\Lambda}_1 $}
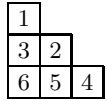
\multiput(0,0)(0,10){2}{\line(1,0){30}}
\put(0,20){\line(1,0){20}}\put(0,30){\line(1,0){10}}
\multiput(0,0)(10,0){2}{\line(0,1){30}}
\put(20,0){\line(0,1){20}}\put(30,0){\line(0,1){10}}
\put(0,0){\makebox(10,10){\scriptsize{$6$}}}
\put(10,0){\makebox(10,10){\scriptsize{$5$}}}
\put(20,0){\makebox(10,10){\scriptsize{$4$}}}
\put(0,10){\makebox(10,10){\scriptsize{$3$}}}
\put(10,10){\makebox(10,10){\scriptsize{$2$}}}
\put(0,20){\makebox(10,10){\scriptsize{$1$}}}
\end{picture}
$$ 
\end{example}
\begin{propos}\label{p:hw}
For each $\Lambda \in P_N^+(L)$ such that $\ov{\Lambda} \equiv \ov{\Lambda}_M\bmod \ov{Q}_N,$ $\psi_{\Lambda}$ is a highest weight vector of $H\otimes \UNp\otimes \ULp.$ The $\UNp$-weight of $\psi_{\Lambda}$ is $\Lambda,$ and  the $\ULp$-weight of $\psi_{\Lambda}$ is $\dot{\Lambda}^{(M)}.$ 
\end{propos}
\begin{proof}
The weights of $\psi_{\Lambda}$ are given by (\ref{eq:wt1}, \ref{eq:wt2}). 
To prove that  $\psi_{\Lambda}$ is annihilated by $\EN_i,\EL_a$ and $B_n$ $(n >0)$ we use the following lemma.
\begin{lemma}\label{l:hw}
Keeping $\Lambda$ as in the statement of Proposition \ref{p:hw}, define the decreasing sequence $k_1,k_2,\dots$ from $\psi_{\Lambda} = u_{k_1}\wedge u_{k_2} \wedge \cdots .$

\noindent Then for  $l>m$ we have
\begin{equation}
u_{k_l} \wedge u_{k_m} = \sum_{l^{\prime }}c_{\alpha , k_{l^{\prime }}}u_{\alpha }\wedge u_{k_{l^{\prime }}} \: \: \mbox{ where } \alpha > k_{l^{\prime }} \geq k_l.
\end{equation}
\end{lemma}
\begin{proof}
Define $\epsilon _{k_i}, a_{k_i}, m_{k_i}$ $(1 \leq \epsilon _{k_i}\leq N, 1\leq  a_{k_i} \leq L, m_{k_i} \in \Zint)$ by  $k_{i} = \epsilon _{k_i}-N(a_{k_i}+L m_{k_i})$.
Using the normal ordering rules, we have 
\begin{equation}
u_{k_l} \wedge u_{k_m} = \sum c_{\alpha ,\beta} u_{\alpha }\wedge u_{\beta}, 
\label{albeta}
\end{equation}
where $k_m\geq \alpha > \beta \geq k_l$ and $\alpha = \epsilon _{k_i} -N(a_{k_{i^{\prime}}}+Lm_{\alpha })$, $\beta = \epsilon _{k_j} -N(a_{k_{j^{\prime}}}+Lm_{\beta })$, $i,j,i^{\prime},j^{\prime} \in \{l,m \} $, $i \neq j$, $i^{\prime} \neq j^{\prime}$, $m_{\alpha }, m_{\beta } \in \Zint$.
\mbox{} From  the explicit expression for  $\psi_{\Lambda}$ (cf. \ref{eq:hwv}) it follows that there is at most one integer $\gamma $ such that $\gamma =  \epsilon _{k_i} -N(a_{k_{i^{\prime}}}+Lm_{\gamma })$ $(i,i^{\prime } \in \{k,l \}, \; m_{\gamma} \in \Zint)$, $k_l < \gamma < k_m $ and $\gamma \neq k_i.$ 
Moreover, if the integer $\gamma $ exists, then $a_l \neq a_m$, $\epsilon _l > \epsilon _m$ and $\gamma = e_{k_l}-N(a_{k_l}+L(m_{k_m}+\delta(a_{k_l} <a_{k_m})))$.
Note that $\gamma$ is the maximal element of the set $\{ \gamma^{\prime} | \gamma^{\prime}  =\epsilon _{k_i} -N(a_{k_{i^{\prime}}}+Lm_{\gamma^{\prime}}), \; i,i^{\prime } \in \{k,l \}, \; m_{\gamma^{\prime}} \in \Zint, \; k_l < \gamma^{\prime} < k_m \} $.
If the $\gamma $ exists, then $\beta$ in (\ref{albeta}) is  distinct from  $\gamma $.
% Because if $\beta = \gamma $ then $\alpha $ must be equal to $k_m$ and we have $\beta =k_l$, it contradicts.
Therefore $\beta = k_{l^{\prime }}$ for some $l^{\prime }$ such that $k_{l^{\prime }} \geq k_l$, and the lemma follows.
\end{proof}
Now we continue the proof of Proposition \ref{p:hw}.
\mbox{} From the definition of $\psi_{\lambda}$ it follows that  $\EN_i\psi_{\Lambda},\EL_a\psi_{\Lambda}$ and $B_n\psi_{\Lambda}$ $(n >0)$ are  linear combinations of vectors of the form 
\begin{equation}
u_{k_1} \wedge \dots\wedge u_{k_{i-1}} \wedge  u_{k_{j}} \wedge  u_{k_{i+1}} \wedge \dots \wedge  u_{k_{j}} \wedge\dots .
\label{klwedge}
\end{equation}
Applying Lemma \ref{l:hw} repeatedly, we conclude that vectors (\ref{klwedge}) are all zero.
\end{proof}

\noindent  Let $\KK[H_-]$ be the Fock module of $H.$ That is $\KK[H_-]$ is the $H$-module generated by the  vector $1$ with the defining relations $B_n 1 =0$ for $n >0.$  

 By Theorem \ref{t:UNUL}, $\F_M$ is an integrable module of $\UNp$ and $\ULp.$ Therefore it is semisimple relative to the algebra $H\otimes \UNp\otimes \ULp.$ Proposition \ref{p:hw} now  implies that  we have an injective $H\otimes \UNp\otimes \ULp$ - linear homomorphism   
\begin{equation}
\bigoplus_{ \{\Lambda \in P_N^+(L)\: |\:  \ov{\Lambda} \equiv \ov{\Lambda}_M\bmod \ov{Q}_N \}} \KK[H_-]\otimes V(\Lambda)\otimes V(\dot{\Lambda}^{(M)}) \; \rightarrow \; \F_M 
\label{eq:hom} \end{equation}
sending $1\otimes v_{\Lambda} \otimes v_{\dot{\Lambda}^{(M)}}$ to $\psi_{\Lambda}.$  It is known (cf. \cite{F1}[Theorem 1.6]) that (\ref{eq:hom}) specializes to an isomorphism when $q=1$. 
The characters of $\KK[H_-],$ $V(\Lambda),$ $V(\dot{\Lambda}^{(M)}),$ and $\F_M$ remain unchanged when $q$ is specialized to $1.$ Therefore (\ref{eq:hom}) is an isomorphism.
Summarizing, we have the following theorem.
\begin{theor} \label{t:decofF}
There is an isomorphism of  $H\otimes \UNp\otimes \ULp$-modules:
\begin{equation}
 \F_M  \cong \bigoplus_{ \{\Lambda \in P_N^+(L)\: |\:  \ov{\Lambda} \equiv \ov{\Lambda}_M\bmod \ov{Q}_N \}} \KK[H_-]\otimes V(\Lambda)\otimes V(\dot{\Lambda}^{(M)}) .
\end{equation}
\end{theor}

\section{The toroidal Hecke algebra and the quantum toroidal algebra} \label{s:tor}
\subsection{Toroidal Hecke algebra} \label{sec:torHecke}
\mbox{} From now on we will work over the base field $\Rat(q^{\frac{1}{2N}})$ rather than $\Rat(q).$ Until the end of the paper we put  $\KK = \Rat(q^{\frac{1}{2N}}).$  Clearly, all results of the preceding sections hold for this $\KK.$   

The toroidal Hecke algebra of type $\gl_n$, $\ddot{\He}_n,$ \cite{VV1,VV2} is a unital associative algebra over $\KK$ with the generators $\xc^{\pm 1},$ $T_i^{\pm 1},X_j^{\pm 1},Y_j^{\pm 1}, 1\leq i < n, 1\leq j \leq n.$ The defining relations involving $T_i^{\pm 1},X_j^{\pm 1}$ are those of the affine Hecke algebra (\ref{eq:ah1} -- \ref{eq:ah2}), and the rest of the relations are as follows:
\begin{eqnarray*}
& \text{ the elements $\xc^{\pm 1}$ are central,}  & \qquad \xc \xc^{-1} = \xc^{-1}\xc = 1,  \\   
& Y_j Y_j^{-1} = Y_j^{-1} Y_j = 1, & \qquad Y_i  Y_j  = Y_j  Y_i, \\   
& T_i^{-1} Y_i T_i^{-1} = q^{-2} Y_{i+1}, &\qquad T_i Y_j = Y_j T_i \quad \text{if $ j \neq i,i+1.$} \\ 
&(X_1 X_2 \cdots X_n)Y_1 = \xc Y_1 (X_1 X_2 \cdots X_n), &\qquad X_2 Y_1^{-1} X_2^{-1} Y_1 = q^{-2} T_1^2.  
\end{eqnarray*}
The subalgebras of  $\ddot{\He}_n$ generated by $T_i^{\pm 1},X_j^{\pm 1}$ and by $T_i^{\pm 1},Y_j^{\pm 1}$ are both isomorphic to the affine Hecke algebra $\dot{\He}_n$ (cf. \cite{VV1}, \cite{VV2}). \\

\noindent Following \cite{C3} we introduce a representation of the toroidal Hecke algebra on the space $ (\KK[z^{\pm 1}] \otimes \KK^L)^{\otimes n}$ $=$ $ \KK[z_1^{\pm 1},\dots,z_n^{\pm 1}] \otimes (\KK^L)^{\otimes n}.$ This representation is  an extension of the representation of $\dot{\He}_n = \langle T_i^{\pm 1} , X_j \rangle $ described in Section \ref{sec:affH}.     

Let $\nu = \sum_{a=1}^L \nu(a) \ep_a,$ where $\ep_a = \dot{\ov{\Lambda}}_a - \dot{\ov{\Lambda}}_{a-1},$ be an integral weight of $\sll_L$ $( \nu(a) \in \Zint).$ Define $q^{\nu^{\vee}} \in \End \left(\KK[z^{\pm 1}] \otimes \KK^L\right)$ as follows:
$$ q^{\nu^{\vee}}( z^m\ef_a) = q^{\nu(L+1-a)} z^m\ef_a.$$   
Here the basis $\ef_1,\dots,\ef_L$ of $\KK^L$ is the same as in Section \ref{sec:Usl}. For $p \in q^{\Zint}$ define $p^D   \in \End \left(\KK[z^{\pm 1}] \otimes \KK^L\right)$ as 
$$ p^{D}( z^m\ef_a) = p^{m} z^m\ef_a.$$   
For $i=1,2,\dots,n-1$ let $s_i$ be the permutation operator of factors $i$ and $i+1$ in   $ (\KK[z^{\pm 1}] \otimes \KK^L)^{\otimes n},$ and let $\tilde{T}_{i,i+1} = - q (\Tc_i)^{-1}.$ Here $\Tc_i$ is the generator of the finite Hecke algebra defined in  (\ref{eq:Tci}). For $X \in \End \left(\KK[z^{\pm 1}] \otimes \KK^L\right)$ let $$(X)_i := 1^{\otimes (i-1)} \otimes X \otimes 1^{\otimes (n-i-1)}\in  \End \left(\KK[z^{\pm 1}] \otimes \KK^L\right)^{\otimes n}.$$ 
For $i=1,2,\dots,n$ define the matrix analogue of the Cherednik-Dunkl operator \cite{C3} as 
\begin{equation}
Y_i^{(n)} = \tilde{T}^{-1}_{i,i+1}\cdots \tilde{T}^{-1}_{n-1,n} s_{n-1} s_{n-2} \cdots s_1 (p^D)_1 (q^{\nu^{\vee}})_1  \tilde{T}_{1,2}\cdots \tilde{T}_{i-1,i}.
\label{eq:Y}
\end{equation}
Let $s\in \{0,1,\dots,NL-1\}$ and $m\in \Zint$ be defined from $n = s + NLm.$ Put $\unun{n} = Nm.$ 
\begin{propos}[\cite{C3}] \label{p:torHeckerep}
The map 
$$ T_i \mapsto \Tc_i,\quad X_i \mapsto z_i, \quad Y_i \mapsto q^{-{\unun{n}}} Y_i^{(n)}, \quad \xc \mapsto p 1 $$ 
extends to a right representation of $\ddot{\He}_n$ on $(\KK[z^{\pm 1}] \otimes \KK^L)^{\otimes n}.$ 
\end{propos}
\noindent{\bf Remark.} The normalizing factor $q^{-{\unun{n}}}$ in the map $Y_i \mapsto q^{-{\unun{n}}} Y_i^{(n)}$ above clearly can be replaced by any coefficient in $\KK.$ The adopted choice of this factor makes $q^{-{\unun{n}}} Y_i^{(n)}$ to behave appropriately (see Proposition \ref{p:inter}) with respect to increments of $n$ by  steps of the value $NL.$   \\

\noindent Let $\chi = \sum_{a=1}^L \chi(a) \ep_a $ be an integral weight of $\sll_L.$ Let  $\UU_q({\mathfrak b}_L)^{\chi}$ be the non-unital subalgebra of $\ULp$ generated by the elements 
\begin{equation}\FL_0,\FL_1,\dots,\FL_{L-1} \quad \text{and} \quad \KL_a - q^{\chi(a) - \chi(a+1)} 1 \quad (a=1,\dots,L-1). 
\label{eq:gen} \end{equation}    
We define an action of $\ULp$ on $\KK[z^{\pm 1}] \otimes \KK^L$ by the obvious restriction of the action on $\KK[z^{\pm 1}] \otimes \KK^L \otimes \KK^N$ defined in (\ref{eq:ul1} -- \ref{eq:ul2}). Iterating the  coproduct $\Delta^-$ given in (\ref{eq:co1} -- \ref{eq:co3}) we obtain an action of  $\ULp$ on $(\KK[z^{\pm 1}] \otimes \KK^L)^{\otimes n}.$ 

\begin{propos} \label{p:inv1}
Suppose $p = q^{-2L},$ and $ \nu = -\chi - 2\rho,$ where $\rho = \sum_{a=1}^{L-1} \dot{\ov{\Lambda}}_a.$ Then the action of the toroidal Hecke algebra on $(\KK[z^{\pm 1}] \otimes \KK^L)^{\otimes n}$ defined in Proposition \ref{p:torHeckerep} leaves invariant the subspace $\UU_q({\mathfrak b}_L)^{\chi} \left((\KK[z^{\pm 1}] \otimes \KK^L)^{\otimes n}\right).$
\end{propos}
\begin{proof}
It is clear that the multiplication by $z_i,$ and hence action of $X_i$ commutes with all generators of  $\ULp.$ \mbox{} From the  intertwining property of the $R$-matrix it follows that the operators $\Tc_i$ (cf. \ref{eq:Tc}) commute with all generators of $\ULp$ as well. With $p = q^{-2L},$ and $ \nu = -\chi - 2\rho,$ a direct computation gives     
\begin{align*}
&Y_n^{(n)} \FL_a = \left( (q^{\chi(a)-\chi(a+1)}1 - \KL_a) \KL_a^{-1} (\FL_a)_n (\KL_a)_n  + \FL_a (\KL_a)_n \right) Y_n^{(n)} \quad (a=1,\dots,L-1),\\
&Y_n^{(n)} \FL_0 = \left( (q^{\chi(L)-\chi(1)}1 - \KL_0) \KL_0^{-1} (\FL_0)_n (\KL_0)_n  + \FL_0 (\KL_0)_n \right) Y_n^{(n)}.
\end{align*}
In view of the relation $\Tc_i Y_{i+1}^{(n)} \Tc_i = q^2 Y_i^{(n)},$ and the commutativity of $\Tc_i$ with the generators of $\ULp,$ this shows that for all $i$ the operators $Y_i^{(n)}$ leave the image of $\UU_q({\mathfrak b}_L)^{\chi} $ invariant.  
\end{proof}

\subsection{The quantum toroidal algebra} \label{sec:tor}
Fix an integer $N\geq 3.$ 
The quantum toroidal algebra of type $\sll_N,$ $\tor,$ is an associative unital algebra over $\KK$ with generators:
$$E_{i,k},\quad F_{i,k},\quad H_{i,l},\quad K_i^{\pm 1}, \quad 
q^{\pm\frac12 c}, \quad \dc^{\pm 1},$$
where  $k\in {\mathbb Z}$, $l\in {\mathbb Z}\backslash \{0\}$ and $i=0,1,\cdots,N-1$. The generators $q^{\pm\frac12 c}$ and $\dc^{\pm 1}$ are central. The rest of the defining relations are expressed in terms of the formal series
$$E_i(z)=\sum_{k\in {\mathbb Z}}E_{i,k}z^{-k}, \quad 
F_i(z)=\sum_{k\in {\mathbb Z}}F_{i,k}z^{-k}, \quad 
K_i^{\pm}(z)=K_i^{\pm 1}\exp(\pm (q-q^{-1})\sum_{k\geq 1}H_{i,\pm k}
z^{\mp k}),$$
as follows:
%\begin{equation}
% q^{\pm\frac12 c}\text{ are central,}
%\label{defrel1}
%\end{equation}
\begin{gather}
 K_i K_i^{-1}=K_i^{-1}K_i= q^{\frac12 c} q^{-\frac12 c} = q^{-\frac12 c} q^{\frac12 c}= \dc \dc^{-1} = \dc^{-1} \dc    = 1, \\
 K_i^{\pm}(z)K_j^{\pm}(w)=K_j^{\pm}(w)K_i^{\pm}(z) 
\label{rb} \\
 \theta_{- a_{ij}}
 (q^{-c}\dc^{m_{ij}}\frac{z}{w})K_i^-(z)K_j^+(w)=
 \theta_{-a_{ij}}
 (q^{c}\dc^{m_{ij}}\frac{z}{w})K_j^+(w)
 K_i^-(z)
\label{rc} \\
 K_i^{\pm}(z)E_j(w)
 =\theta_{\mp a_{ij}}
 (q^{-\frac12 c}\dc^{\mp m_{ij}}w^{\pm}z^{\mp})E_j(w)K_i^+(z)
\label{rd} \\
 K_i^{\pm}(z)F_j(w)
 =\theta_{\pm a_{ij}}
 (q^{\frac12 c}\dc^{\mp m_{ij}}w^{\pm}z^{\mp})F_j(w)K_i^+(z) \\
 [E_i(z),F_j(w)]=\delta_{i,j}\frac{1}{q-q^{-1}}
 \{\delta(q^c\frac{w}{z})K_i^+(q^{\frac12 c}w)-\delta(q^c\frac{z}{w})K_i^
 -(q^{\frac12 c}z)\}
\label{re} \\
 (\dc^{m_{ij}}z-q^{a_{ij}}w)E_i(z)E_j(w)
 =(q^{a_{ij}}\dc^{m_{ij}}z-w)E_j(w)E_i(z)
\label{rf} \\
 (\dc^{m_{ij}}z-q^{-a_{ij}}w)F_i(z)F_j(w)
 =(q^{-a_{ij}}\dc^{m_{ij}}z-w)F_j(w)F_i(z) \\
 \sum_{\sigma\in {\mathfrak S}_m}\sum_{r=0}^m(-1)^r
 \begin{bmatrix}m\\r\end{bmatrix}
 E_i(z_{\sigma(1)})\cdots E_i(z_{\sigma(r)})E_j(w)E_i(z_{\sigma(r+1)})\cdots
 E_i(z_{\sigma(m)})=0
\label{rg} \\
  \sum_{\sigma\in {\mathfrak S}_m}\sum_{r=0}^m(-1)^r
 \begin{bmatrix}m\\r\end{bmatrix}
{} F_i(z_{\sigma(1)})\cdots F_i(z_{\sigma(r)})F_j(w)F_i(z_{\sigma(r+1)})\cdots
{} F_i(z_{\sigma(m)}) = 0 \label{rg1}
\end{gather}
where in (\ref{rg}) and (\ref{rg1})   $i\ne j$ and $m=1-a_{ij}$.\\

\noindent In these  defining relations $\delta(z) = \sum_{n = -\infty}^{\infty} z^n,$ $\theta_m(z)  \in \KK[[z]] $ is the expansion of $\frac{zq^m-1}{z-q^m},$ $a_{ij}$ are the entries of the Cartan matrix of $\asll_N,$ and  
$m_{ij}$ are the entries
of the following $N\times N$-matrix
$$
M=\begin{pmatrix}
      0 &     -1 &      0 & \hdots &      0 &      1\\
      1 &      0 &     -1 & \hdots &      0 &      0\\
      0 &      1 &      0 & \hdots &      0 &      0\\
 \vdots & \vdots & \vdots & \ddots & \vdots & \vdots\\
      0 &      0 &      0 & \hdots &      0 &     -1\\
     -1 &      0 &      0 & \hdots &      1 &      0
\end{pmatrix}.
$$
Let $\UU_h$ be the subalgebra of $\tor$ generated by the elements $E_{i,0},F_{i,0},K_i^{\pm 1}$ $(0\leq i < N).$ These elements satisfy the defining relations (\ref{eq:r1} -- \ref{eq:r3}) and (\ref{eq:r5} -- \ref{eq:r7})  of $\UNp.$ Thus the following map extends to a homomorphism of algebras: 
\begin{equation}\UNp \rightarrow \UU_h :\; \EN_i \mapsto  E_{i,0}, \quad \FN_i \mapsto  F_{i,0}, \quad \KN_i^{\pm 1} \mapsto  \KN_i^{\pm 1}. \label{eq:Uh}
\end{equation}
Let  $\UU_v$ be the subalgebra of $\tor$ generated by the elements $E_{i,k},F_{i,k},H_{i,l},$ $K_i^{\pm 1}$ $(1\leq i < N; k\in \Zint; l\in \Zint_{\neq 0}),$ and $q^{\pm \frac12 c}, \dc^{\pm 1}.$ Recall, that apart from the presentation given in Section \ref{sec:Usl}, the algebra $\UNp$ has the ``new presentation'' due to Drinfeld which is similar to that one of $\tor$ above. A proof of the  isomorphism between the two presentations is announced in  \cite{Drinfeld1} and given in \cite{Beck}. Let $\Et_{i,k},\Ft_{i,k},\Ht_{i,l},\Kt_{i}^{\pm 1},$ $(1\leq i < N; k\in \Zint; l\in \Zint_{\neq 0}),$ and $ q^{\pm \frac12 \ct}$ be the generators of $\UNp$ in the realization of \cite{Drinfeld1}. Comparing this realization of  $\UNp$ with the defining relations of $\tor$ one easily sees that the map      
\begin{align}
\UNp \rightarrow \UU_v :\; &\Et_{i,k}  \mapsto  \dc^{ik}E_{i,k},\quad \Ft_{i,k}  \mapsto  \dc^{ik}F_{i,k},\quad \Ht_{i,l}  \mapsto  \dc^{il}H_{i,l},\label{eq:Uv}\\ &\Kt_i^{\pm 1} \mapsto  \KN_i^{\pm 1}, \quad  q^{\pm \frac12 \ct} \mapsto q^{\pm \frac12 c}  \nonumber 
\end{align}
where $1 \leq i < N,$ extends to a homomorphism of algebras. 
Thus each module of $\tor$ carries two actions of $\UNp$ obtained by pull-backs through the homomorphisms (\ref{eq:Uh}) and (\ref{eq:Uv}). We will say that a module of $\tor$ has {\em level} $(l_v,l_h)$ provided the action of $\UNp$ obtained through the homomorphism (\ref{eq:Uh}) has level $l_h,$ and the action of $\UNp$ obtained through the homomorphism (\ref{eq:Uv}) has level $l_v.$ On such a module the central elements $q^{\pm \frac12 c}$ act as multiplications by $q^{\pm \frac12 l_v},$ and the element $K_0 K_1 \cdots K_{N-1}$ acts as the multiplication by $q^{l_h}.$   

The following proposition, proved in \cite{VV1}, shows that it is sometimes possible to extend  a representation of $\UNp$ to a representation of $\tor.$ 
\begin{propos} \label{p:shift1}
Let $W$ be a module of $\UNp.$ Suppose that there are $a,b \in q^{\Zint},$ and an invertible $\tilde{\psi} \in \End(W)$ such that  
\begin{alignat}{5}
&\tilde{\psi}^{-1}\Et_i(z)\tilde{\psi} = \Et_{i-1}(az),& &\tilde{\psi}^{-2}\Et_1(z)\tilde{\psi}^2 = \Et_{N-1}(bz),&  \\   
&\tilde{\psi}^{-1}\Ft_i(z)\tilde{\psi} = \Ft_{i-1}(az),& &\tilde{\psi}^{-2}\Ft_1(z)\tilde{\psi}^2 = \Ft_{N-1}(bz),& \\
&\tilde{\psi}^{-1}\Kt_i^{\pm}(z)\tilde{\psi} = \Kt_{i-1}^{\pm}(az),& \qquad&\tilde{\psi}^{-2}\Kt_1^{\pm}(z)\tilde{\psi}^2 = \Kt_{N-1}^{\pm}(bz),&
\end{alignat}
where $2\leq i < N.$ Then $W$ is a $\tor$-module with the action given by 
\begin{align*}
& X_i(z) = \Xt_i(d^iz) \qquad  (1\leq i <N), \qquad   
 X_0(z) = \tilde{\psi}^{-1}\Xt_1(a^{-1}d^{-1}z)\tilde{\psi}, \\
& \dc  = d 1, \qquad q^{\frac12 c} = q^{\frac12 \ct}.  
\end{align*}
where $d^N = b/a^2,$ and $X = E,F,K^{\pm}.$
\end{propos}

\subsection{The Varagnolo-Vasserot duality} \label{sec:VVdual} We now briefly review, following \cite{VV1}, the Schur-type duality between the toroidal Hecke algebra $\ddot{\He}_n$ and the quantum toroidal algebra $\tor.$ 

Let $\M$ be a right $\ddot{\He}_n $-module, such  that the central element $\xc$ of $\ddot{\He}_n$ acts as the multiplication by $x \in q^{\Zint}.$ The algebra $\ddot{\He}_n$ contains two subalgebras: $\dot{\He}_n^h = \langle T_i^{\pm 1} , X_j \rangle ,$ and $\dot{\He}_n^v = \langle T_i^{\pm 1} , Y_j \rangle$ both isomorphic to the affine Hecke algebra $\dot{\He}_n.$ Therefore the duality functor of Chari--Pressley \cite{CP} yields two actions of $\UNp$ on the linear space 
$ \M\otimes_{\He_n} (\KK^N)^{\otimes n}.$ Here the action of the finite Hecke algebra $\He_n$ on $(\KK^N)^{\otimes n}$ is given by (\ref{eq:Tsi}), and $\He_n$ is embedded into $\ddot{\He}_n$ as the subalgebra generated by $T_i^{\pm 1}.$ 

For $i,j=1,\dots,N$ let $e_{i,j} \in \End(\KK^N)$ be the matrix units with respect to the basis $\vf_1,\vf_2,\dots,\vf_{N}$ (cf. Section \ref{sec:Usl}). For $i=0,1,\dots,N-1$ let $k_i = q^{e_{i,i} - e_{i+1,i+1}},$ where the indices are cyclically extended modulo $N$. For $X \in \End(\KK^N)$ we put $(X)_i = 1^{\otimes (i-1)}\otimes X \otimes 1^{\otimes (n-i)}.$ 

The functor of \cite{CP} applied to $\M$ considered as the $\dot{\He}_n^h$-module gives the following action of $\UNp$ on  $ \M\otimes_{\He_n} (\KK^N)^{\otimes n}:$ 
\begin{align}      
&E_i(m\otimes v) = \sum_{j=1}^n m X_j^{\delta(i=0)} \otimes (e_{i,i+1})_j (k_i)_{j+1} (k_i)_{j+2}\cdots (k_i)_{n}v, \label{eq:h1}\\ 
&F_i(m\otimes v) = \sum_{j=1}^n m X_j^{-\delta(i=0)} \otimes (e_{i+1,i})_j (k_i^{-1})_{1} (k_i^{-1})_{2}\cdots (k_i^{-1})_{j-1}v, \\
&K_i(m\otimes v) = m \otimes (k_i)_1(k_i)_2 \cdots (k_i)_n v. \label{eq:h3}
\end{align}
Here $m \in \M, v \in (\KK^N)^{\otimes n},$ and the indices are cyclically extended modulo $N.$ Likewise, application of this functor to  $\M$ considered as the $\dot{\He}_n^v$-module gives another action of $\UNp$ on  $ \M\otimes_{\He_n} (\KK^N)^{\otimes n}:$
\begin{align}      
&\Ed_i(m\otimes v) = \sum_{j=1}^n m Y_j^{-\delta(i=0)} \otimes (e_{i,i+1})_j (k_i)_{j+1} (k_i)_{j+2}\cdots (k_i)_{n}v,  \label{eq:v1}\\ 
&\Fd_i(m\otimes v) = \sum_{j=1}^n m Y_j^{\delta(i=0)} \otimes (e_{i+1,i})_j (k_i^{-1})_{1} (k_i^{-1})_{2}\cdots (k_i^{-1})_{j-1}v, \\
&\Kd_i(m\otimes v) = m \otimes (k_i)_1(k_i)_2 \cdots (k_i)_n v.  \label{eq:v3}
\end{align}
Here we put hats over the generators in order to distinguish the actions given by (\ref{eq:h1} -- \ref{eq:h3}) and (\ref{eq:v1} -- \ref{eq:v3}).  

Varagnolo and Vasserot have proven, in \cite{VV1}, that $ \M\otimes_{\He_n} (\KK^N)^{\otimes n}$ is a $\tor$-module such that the $\UNp$-action (\ref{eq:h1} -- \ref{eq:h3}) is the pull-back through the homomorphism (\ref{eq:Uh}), and the $\UNp$-action (\ref{eq:v1} -- \ref{eq:v3}) is the pull-back through the homomorphism (\ref{eq:Uv}). Let us recall here the main element of their proof.  

Let $\psi$ be the endomorphism of $ \M\otimes_{\He_n} (\KK^N)^{\otimes n}$ defined by 
\begin{gather}
\psi : m \otimes \vf_{\ep_1}\otimes \vf_{\ep_2}\otimes \cdots \otimes \vf_{\ep_n} \mapsto\label{eq:psi} \\ m X_1^{-\delta_{N,\ep_1}}X_2^{-\delta_{N,\ep_2}} \cdots  X_n^{-\delta_{N,\ep_n}}  \otimes \vf_{\ep_1+1}\otimes \vf_{\ep_2+1}\otimes \cdots \otimes \vf_{\ep_n+1}, \nonumber 
\end{gather}
where $\vf_{N+1}$ is identified with $\vf_1.$ Taking into account the defining relations of $\Htor$ one can confirm that $\psi$ is well-defined.   

Let $\Et_{i,k},\Ft_{i,k},\Ht_{i,l},\Kt_{i}^{\pm 1}$ $(k\in \Zint;l\in \Zint_{\neq 0}; 1\leq i < N)$ be the generators of the $\UNp$-action (\ref{eq:v1} -- \ref{eq:v3}) obtained from $\Ed_j,\Fd_j,\Kd_j^{\pm 1}$ $(0\leq j < N)$ by the isomorphism between the two realizations of $\UNp$ given in \cite{Beck}. Let $\Et_i(z), \Ft_i(z), \Kt^{\pm}_i(z)$ be the corresponding generating series.  

\begin{propos}[\cite{VV1}] \label{p:fintwist}
The following relations hold in $\M\otimes_{\He_n} (\KK^N)^{\otimes n}:$ 
\begin{alignat}{5}
&{\psi}^{-1}\Et_i(z){\psi} = \Et_{i-1}(q^{-1}z),& &{\psi}^{-2}\Et_1(z){\psi}^2 = \Et_{N-1}(x^{-1}q^{N-2}z),&  \\   
&{\psi}^{-1}\Ft_i(z){\psi} = \Ft_{i-1}(q^{-1}z),& &{\psi}^{-2}\Ft_1(z){\psi}^2 = \Ft_{N-1}(x^{-1}q^{N-2}z),& \\
&{\psi}^{-1}\Kt_i^{\pm}(z){\psi} = \Kt_{i-1}^{\pm}(q^{-1}z),& \qquad&{\psi}^{-2}\Kt_1^{\pm}(z){\psi}^2 = \Kt_{N-1}^{\pm}(x^{-1}q^{N-2}z).&
\end{alignat}
Here $2\leq i < N.$ 
\end{propos}
\noindent Proposition \ref{p:shift1} now implies that  $\M\otimes_{\He_n} (\KK^N)^{\otimes n}$ is a $\tor$-module, in particular, the central element $\dc$ acts as the multiplication by $x^{-1/N}q,$ and the central element $q^{\frac12 c}$ acts as the multiplication by $1.$  

\subsection{The action of the quantum toroidal algebra on the wedge product}
In the framework of the preceding section, let $\M = (\KK[z^{\pm 1}]\otimes \KK^L)^{\otimes n}$ be the $\Htor$-module with the action given in Proposition \ref{p:torHeckerep}. In view of the remark made in Section \ref{sec:wedge}, the linear space $\M\otimes_{\He_n} (\KK^N)^{\otimes n}$ is isomorphic to the wedge product $\wedge^n\Vaff.$ Therefore, by the Varagnolo-Vasserot duality,  $\wedge^n\Vaff$ is a module of $\tor.$ The action of $\UNp$ given by (\ref{eq:h1} -- \ref{eq:h3}) coincides with the action of $\UNp$ defined on $\wedge^n\Vaff$ in Section \ref{sec:wedge}. Following the terminology of \cite{VV2}, we will call this action {\em the horizontal} action of  $\UNp$ on $\wedge^n \Vaff.$ The formulas (\ref{eq:v1} -- \ref{eq:v3}) give another action of $\UNp$ on $\wedge^n \Vaff,$ we will refer to this action as {\em the vertical } action. 

Recall, that in Section \ref{sec:wedge} an action of $\ULp,$ commutative with the horizontal action of $\UNp,$ was defined on $\wedge^n\Vaff.$ Recall, as well, that for each integral weight $\chi$ of $\sll_L$ we have defined, in Section \ref{sec:torHecke}, the subalgebra $\UU_q({\mathfrak b}_L)^{\chi} $ of $\ULp.$    
The $\Htor$-module structure defined in Proposition \ref{p:torHeckerep} depends on two parameters: $\nu$ which is  an integral weight of $\sll_L,$ and $p \in q^{\Zint}.$ The same parameters thus enter into the $\tor$-module structure on $\wedge^n \Vaff.$  
\begin{propos}\label{p:inv2}
Suppose $p = q^{-2L},$ and $\nu = - \chi - 2\rho$ for an integral $\sll_L$-weight $\chi.$ Then the action of $\tor$ on $\wedge^n \Vaff$ leaves invariant the linear subspace $\UU_q({\mathfrak b}_L)^{\chi}\left(\wedge^n \Vaff \right).$   
\end{propos}
\begin{proof}
It is not difficult to see, that the subalgebras $\UU_h$ and $\UU_v$ generate $\tor$ (cf. Lemma 2  in \cite{STU}). Therefore, to prove the proposition, it is enough to show, that both the horizontal and the vertical actions of $\UNp$ on $\wedge^n \Vaff$ leave $\UU_q({\mathfrak b}_L)^{\chi}\left(\wedge^n \Vaff \right)$ invariant. However, the horizontal action commutes with the action of $\ULp,$ while Proposition \ref{p:inv1} implies that the vertical action leaves $\UU_q({\mathfrak b}_L)^{\chi}\left(\wedge^n \Vaff \right)$ invariant.    
\end{proof}

\section{The actions of the quantum toroidal algebra on  the Fock spaces and on irreducible integrable highest weight modules of $\UU_q^{\prime}(\agl_N)$} \label{s:toract}

\subsection{A level 0 action of $\UNp$ on the Fock space}
Let $\pi^v_{(n)} : \UNp \rightarrow \End(\wedge^n \Vaff)$ be the map defining the vertical action of $\UNp$ on the wedge product $\wedge^n\Vaff.$ In accordance with (\ref{eq:v1} -- \ref{eq:v3}), for $f\in (\KK[z^{\pm 1}]\otimes \KK^L)^{\otimes n}$ and $v \in (\KK^N)^{\otimes n }$ we have  
\begin{align}      
&\pi^v_{(n)}(\EN_i)\cdot \wedge(f\otimes v) = \wedge \sum_{j=1}^n  (q^{-\unun{n}}Y^{(n)}_j)^{-\delta(i=0)} f \otimes (e_{i,i+1})_j (k_i)_{j+1} (k_i)_{j+2}\cdots (k_i)_{n}v,  \label{eq:pi1}\\ 
&\pi^v_{(n)}(\FN_i)\cdot \wedge(f\otimes v) = \wedge \sum_{j=1}^n  (q^{-\unun{n}}Y^{(n)}_j)^{\delta(i=0)} f \otimes (e_{i+1,i})_j (k_i^{-1})_{1} (k_i^{-1})_{2}\cdots (k_i^{-1})_{j-1}v,\\
&\pi^v_{(n)}(\KN_i)\cdot \wedge(f\otimes v) = \wedge f \otimes (k_i)_1(k_i)_2 \cdots (k_i)_n v,  \label{eq:pi3}
\end{align}
where we denote by $\wedge$ the canonical map from $\Vaff^{\otimes n}$ $ = $ $ 
(\KK[z^{\pm 1}]\otimes \KK^L)^{\otimes n} \otimes (\KK^N)^{\otimes n }$ to $\wedge^n \Vaff.$ 

In this section, for each $M\in \Zint,$ we define a level 0 action of $\UNp$ on the Fock space $\F_M.$ Informally, this action arises as the limit $n\rightarrow \infty$ of the vertical action (\ref{eq:pi1} -- \ref{eq:pi3}) on the wedge product. In parallel with the finite case, the Fock space, thus admits two actions of $\UNp:$ the level $L$ action defined in Section \ref{sec:UF} as the inductive limit of the horizontal action, and an extra action with level zero.     

We start by introducing a grading on $\F_M.$ To facilitate this, we adopt the following notational convention. For each integer $k$ we define the unique triple $\ov{k},\dot{k},\un{k},$ where $\ov{k} \in \{1,2,\dots ,N\},$ $\dot{k} \in \{1,2,\dots ,L\},$ $\un{k} \in \Zint$ by 
$$ k = \ov{k} - N( \dot{k} + L \un{k}) .$$
Then (cf. Section \ref{sec:nowedges}) we have $u_k = z^{\un{k}}\ef_{\dot{k}}\vf_{\ov{k}}.$ The Fock space $\F_M$ has a basis formed by normally ordered semi-infinite wedges  
$ u_{k_1}\wedge u_{k_2} \wedge \cdots $ 
where the decreasing sequence of momenta $k_1,k_2,\dots$ satisfies the asymptotic condition $k_i = M-i+1$ for $i\gg 1.$ Let $o_1,o_2,\dots$ be the sequence of momenta labeling the vacuum vector $|M\rangle$ of $\F_M,$ i.e.: $o_i = M-i+1$ for all $i \geq 1.$ Define the degree of a semi-infinite normally ordered wedge by    
\begin{equation}
\deg u_{k_1}\wedge u_{k_2} \wedge \cdots = \sum_{i\geq 1} \un{o_i} - \un{k_i}.  
\label{eq:sideg}\end{equation}
Let $\F_M^d$ be the homogeneous component of $\F_M$ of degree $d.$ Clearly, the asymptotic condition $k_i = M-i+1$ $(i\gg 1)$ implies that  
$$ \F_M = \bigoplus_{d=0}^{\infty}  \F_M^d.$$ 
Let $s \in \{0,1,\dots,NL-1\}$ be defined from $M\equiv s\bmod NL.$ For a non-negative integer $l$ we define the linear subspace $V_{M,s+lNL}$ of $\wedge^{s+lNL}\Vaff$ by   
\begin{equation}
V_{M,s+lNL} = \bigoplus_{\un{k_{s+lNL}} \leq \un{o_{s+lNL}}} \KK u_{k_1}\wedge u_{k_2}\wedge \cdots \wedge u_{k_{s+lNL}}, \label{eq:VM}
\end{equation}
where the wedges in the right-hand side are assumed to be normally ordered. For $s=l=0$ we put $V_{M,s+lNL} = \KK.$ The vector space (\ref{eq:VM}) has a grading similar to that one of the Fock space. Now the degree of a normally ordered wedge is defined as  
\begin{equation}
\deg u_{k_1}\wedge u_{k_2} \wedge \cdots\wedge u_{k_{s+lNL}} = \sum_{i=1}^{s+lNL} \un{o_i} - \un{k_i}.  
\label{eq:fdeg}\end{equation}
Note that this degree is necessarily a non-negative integer since $k_1 > k_2 > \cdots > k_{s+lNL}$ and $\un{k_{s+lNL}} \leq \un{o_{s+lNL}}$ imply $\un{k_i} \leq \un{o_i}$ for all $i=1,2,\dots,s+lNl.$ Let $V_{M,s+lNL}^d$ be the homogeneous component of $V_{M,s+lNL}$ of degree $d.$  

For non-negative integers  $d$ and $l$ introduce the following linear map:
\begin{equation}
\vro_l^d : V_{M,s+lNL}^d \rightarrow \F_M^d \: : \: w \mapsto w\wedge | M-s-lNL \rangle. 
\label{eq:rol}\end{equation}
The proof of the following proposition is straightforward (cf. Proposition 16 in \cite{STU}, or Proposition 3.3 in \cite{U}). 
\begin{propos} \label{p:rho}
Suppose $l \geq d.$ Then $\vro_l^d$ is an isomorphism of vector spaces. 
\end{propos}
\noindent In view of  this proposition, it is clear that for non-negative integers $d,l,m,$ such that  $d \leq l < m,$ the linear map   
\begin{equation}
\vro_{l,m}^d : V_{M,s+lNL}^d \rightarrow V_{M,s+mNL}^d \: : \: w \mapsto w\wedge u_{M-s-lNL}\wedge u_{M-s-lNL-1}\wedge \cdots \wedge u_{M-s-mNL+1}   
\label{eq:rolm}\end{equation}
is an isomorphism of vector spaces as well.

Now let us return to the vertical action $\pi^v_{(n)}$ of $\UNp$ on  $\wedge^n \Vaff $ given by (\ref{eq:pi1} -- \ref{eq:pi3}). 
\begin{propos}
For each $d=0,1,\dots$ the subspace $ V_{M,s+lNL}^d  \subset \wedge^{s+lNL} \Vaff $ is invariant with respect to the action $\pi^v_{(s+lNL)}.$ 
\end{propos} 
\begin{proof}
Let $n=s+lNL,$ and let us identify $\Vaff^{\otimes n}$ with $\KK[z_1^{\pm 1},\dots,z_n^{\pm 1}]\otimes (\KK^L)^{\otimes n} \otimes (\KK^N)^{\otimes n}$ by the  isomorphism 
$$ z^{m_1}\ef_{a_1}\vf_{\ep_1} \otimes \cdots \otimes z^{m_n}\ef_{a_n}\vf_{\ep_n} \mapsto z_1^{m_1}\cdots z_n^{m_n} \ef_{a_1} \cdots \ef_{a_n} \vf_{\ep_n} \cdots \vf_{\ep_n}. $$ 
Then $V_{M,s+lNL}$ is the image, with respect to the quotient map $\wedge : \Vaff^{\otimes n} \rightarrow \wedge^n \Vaff,$  of the subspace
\begin{equation} (z_1 \cdots z_n)^{\un{o_n}} \KK[z_1^{-1},\dots,z_n^{- 1}]\otimes (\KK^L)^{\otimes n} \otimes (\KK^N)^{\otimes n} \subset \Vaff^{\otimes n},  
\label{eq:ss}\end{equation}
while the grading on $V_{M,s+lNL}$ is induced from the grading of (\ref{eq:ss}) by eigenvalues of the operator $D = z_1\frac{\partial}{\partial z_1} + \cdots + z_n\frac{\partial}{\partial z_n}.$ 

The operators $Y_i^{(n)}$ leave $(z_1 \cdots z_n)^{\un{o_n}} \KK[z_1^{-1},\dots,z_n^{- 1}]\otimes (\KK^L)^{\otimes n} $ invariant, and commute with $D.$ Now (\ref{eq:pi1} -- \ref{eq:pi3}) imply the statement of the proposition.     
\end{proof}

\begin{propos} \label{p:inter}
Let $0 \leq d \leq l,$ let $n=s+lNL,$ and let $X$ be any of the generators $E_i,F_i,K_i^{\pm 1}$ $(0\leq i < N)$ of $\UNp.$ Then the following intertwining relation holds for all  $w \in V_{M,s+lNL}^d:$
\begin{equation}
\pi^v_{(n+NL)}(X) \cdot \vro_{l,l+1}^d (w) = \vro_{l,l+1}^d \left(\pi^v_{(n)}(X)\cdot w \right).
\label{eq:inter}\end{equation}
Consequently, for $0 \leq d \leq l < m$ the map $\vro_{l,m}^d$ defined in {\em (\ref{eq:rolm})} is an isomorphism of $\UNp$-modules. 
\end{propos}
\begin{proof}
The proof is based, in particular, on Lemma \ref{l:ml}, to state which we introduce the following notation. For $\mm = (m_1,m_2,\dots,m_n) \in \Zint^n,$ and $\aab = (a_1,a_2,\dots,a_n) \in \{1,2,\dots,L\}^n$ let 
$$ \zeta_i(\mm,\aab) = p^{m_i} q^{\nu(L+1-a_i) + \mu_i(\mm,\aab)} \qquad (i=1,2,\dots,n) $$ 
where $p, \nu$ are the parameters of the representation of $\Htor$ introduced in Section \ref{sec:torHecke}, and $\mu_i(\mm,\aab) = - \#\{j < i | m_j < m_i, a_j = a_i\} + \#\{j < i | m_j \geq  m_i, a_j = a_i\} + \#\{j > i | m_j > m_i, a_j = a_i\} - \#\{j > i | m_j \leq  m_i, a_j = a_i\}.$ 

\begin{lemma}\label{l:ml}
For $k =1,2,\dots$ consider the following monomial 
$$ f =  z_1^{m_1} z_2^{m_2} \dots z_{n+k}^{m_{n+k}} \otimes \ef_{a_1}\ef_{a_2} \cdots \ef_{a_{n+k}} \in \KK[z_1^{\pm 1},\dots,z_{n+k}^{\pm 1}]\otimes (\KK^L)^{\otimes(n+k)}.$$ 
Assume that $m_1,m_2,\dots,m_n < m_{n+1} = m_{n+2} = \cdots = m_{n+k} =: m,$ and that $a_{n+i} \leq a_{n+j}$ for $ 1 \leq i < j \leq k.$ For $j \in \{ 1,2,\dots,L\}$ put $\ov{n}(j) = \#\{ \: i\: | \: a_{n+i} = j, 1\leq i \leq k\}.$  

Define the linear  subspaces  $\KC_{n,k}^m,\LC_{n,k}^m  \subset  \KK[z_1^{\pm 1},\dots,z_{n+k}^{\pm 1}]\otimes (\KK^L)^{\otimes(n+k)}$ as follows:
\begin{align*}
&\KC_{n,k}^m = \KK\{ z_1^{m_1^{\prime}}\cdots z_{n+k}^{m_{n+k}^{\prime}} \otimes \ef \: | \: \ef \in (\KK^L)^{n+k}; m_1^{\prime},\dots,m_{n+k}^{\prime} \leq m; \#\{ m_i^{\prime} | m_i^{\prime} = m \} < k\}, \\    
&\LC_{n,k}^m = \KK\{ z_1^{m_1^{\prime}}\cdots z_{n+k}^{m_{n+k}^{\prime}} \otimes \ef_{b_1}\cdots \ef_{b_{n+k}} \: |\: m_1^{\prime},\dots,m_{n}^{\prime} < m ;  m_{n+1}^{\prime},\dots,m_{n+k}^{\prime} = m; \\ & \qquad \qquad \qquad \qquad \qquad \qquad \qquad \qquad \exists j < a_{n+k} \:\text{{\em s.t.}}\: \#\{\: i \: | b_{n+i} = j, 1 \leq i \leq k\} > \ov{n}(j)\}.  
\end{align*} 
Then 
\begin{align*}
&\left( Y_i^{(n+k)} \right)^{\pm 1}(f) \equiv \zeta_i(\mm,\aab)^{\pm 1} f \:\bmod \left( \KC_{n,k}^m  + \LC_{n,k}^m \right) \quad (i=n+1,n+2,\dots,n+k), \\ 
&\left( Y_i^{(n+k)} \right)^{\pm 1}(f) \equiv  q^{\pm \ov{n}(a_i)} \left( Y_i^{(n)} \right)^{\pm 1}(f) \: \bmod \left( \KC_{n,k}^m  + \LC_{n,k}^m \right) \quad (i=1,2,\dots,n). 
\end{align*}
Here $\mm = (m_1,\dots,m_{n+k}),$  $\aab = (a_1,\dots,a_{n+k}),$ and in the right-hand side of the last equation  $\left( Y_i^{(n)} \right)^{\pm 1}$ act on the first $n$ factors of the monomial $f.$ 
\end{lemma}
\noindent A proof of the lemma is given in \cite{TU} for $L=1.$ A proof for general $L$ is quite similar and will be omitted here.  

Let $w$ be a normally ordered wedge from $V_{M,n}^d,$ and let $\bar{w} = \vro_{l,l+1}^d(w).$ The vector $\bar{w}$ is a normally ordered wedge from $V_{M,n+NL}^d,$ we have 
\begin{equation}
\bar{w} = u_{k_1}\wedge u_{k_2} \wedge \cdots \wedge u_{k_{n+NL}} = \wedge (f \otimes v),
\end{equation}
where 
\begin{align}
& f = (z_1^{\un{k_1}}z_2^{\un{k_2}} \cdots z_n^{\un{k_n}})(z_{n+1} \cdots z_{n+NL})^{m}\otimes \label{eq:fff}\\ & \qquad \qquad \qquad \qquad\otimes (\ef_{\dot{k}_1}\ef_{\dot{k}_2} \cdots \ef_{\dot{k}_n})\underbrace{(\ef_1 \cdots \ef_1)}_{ N \: {\mathrm {times}}}\underbrace{(\ef_2 \cdots \ef_2)}_{ N \:{\mathrm {times}}} \dots \underbrace{(\ef_L \cdots \ef_L)}_{ N \:{\mathrm {times}}}, \nonumber \\  
& v = (\vf_{\ov{k_1}}\vf_{\ov{k_2}} \cdots \vf_{\ov{k_n}})(\vf_N \vf_{N-1} \underbrace{ \cdots \vf_1 ) \dots (\vf_N \vf_{N-1} }_{ L \: {\mathrm {copies}}} \cdots \vf_1) \in (\KK^N)^{\otimes (n+NL)}, 
\end{align}
and $m = \un{o_{n+1}} =  \un{o_{n+2}} = \cdots =\un{o_{n+NL}}.$ The monomial $f$ given by (\ref{eq:fff}) satisfies the assumptions of Lemma \ref{l:ml} with $k=NL,$ and $\ov{n}(j) = N$ for all $j \in \{1,2,\dots,L\}.$ Let 
$\KC_{n,NL}^m$ and $\LC_{n,NL}^m$ be the corresponding subspaces of $\KK[z_1^{\pm 1},\dots,z_{n+NL}^{\pm 1}]\otimes (\KK^L)^{\otimes(n+NL)}.$
\begin{lemma} \label{l:ml1}
Let $y \in (\KK^N)^{\otimes (n+NL)},$ and let $f_1 \in \KC_{n,NL}^m,$ $f_2 \in \LC_{n,NL}^m. $ Then 
\begin{align}
 &\wedge( f_1 \otimes y ) \in \oplus_{d^{\prime} > l} V_{M,n+NL}^{d^{\prime}}, \tag{i} \\ 
 &\wedge( f_2 \otimes y )  = 0. \tag{ii}
\end{align}
\end{lemma}
\begin{proof} \mbox{} \\ 
This lemma is the special case ($b=L$ and $c=N$) of Lemma \ref{l:mlbc}. See the proof of Lemma \ref{l:mlbc}.
\end{proof}
Now we continue the proof of the proposition. \mbox{} From the definitions (\ref{eq:pi1} -- \ref{eq:pi3}) and  Lemmas \ref{l:lemma}, \ref{l:ml} and \ref{l:ml1}, it follows that (\ref{eq:inter}) holds modulo $\oplus_{d^{\prime} > d} V_{M,n+NL}^{d^{\prime}}.$ However, the both sides of (\ref{eq:inter}) belong to $V_{M,n+NL}^{d}$ since the action of $\UNp$ preserves the degree $d.$ Hence (\ref{eq:inter}) holds exactly.    

\end{proof} 
Now we are ready to give the definition of the level 0 action of $\UNp$ on the Fock space $\F_M.$ 

\begin{defin} \label{dp:def}
Let $0 \leq d \leq l.$ We define a $\UNp$-action $\pi^v : \UNp \mapsto \End(\F_M^d)$ as 
$$ \pi^v(X) = \vro_l^d \circ \pi^v_{(s+lNL)}(X) \circ (\vro_l^d)^{-1} \qquad (X \in \UNp ).$$
By Proposition  \ref{p:inter} this definition does not depend on the choice of $l$ as long as $l \geq d.$ 
\end{defin}
Thus a $\UNp$-action is defined on each homogeneous component $\F_M^d,$ and hence on the entire Fock space $\F_M.$ 

\subsection{The action of the quantum toroidal algebra on the Fock space}
In Section \ref{sec:UF} we defined a level $L$ action of $\UNp$ on $\F_M.$ Let us denote by $\pi^h$ the corresponding map $\UNp \rightarrow \End(\F_M).$ We refer to  $\pi^h$ as the horizontal $\UNp$-action on the Fock space. In the preceding section we defined another -- level 0 -- action  $ \pi^v  : \UNp \rightarrow \End(\F_M).$ We call $\pi^v$ the vertical $\UNp$-action. Note that for $i=1,2,\dots,N-1$ we have 
\begin{equation*}
\pi^h(\EN_i) = \pi^v(\EN_i), \quad  
\pi^h(\FN_i) = \pi^v(\FN_i), \quad 
\pi^h(\KN_i) = \pi^v(\KN_i),
\end{equation*}
i.e. the restrictions of $\pi^h$ and $\pi^v$ on the subalgebra $\UU_q(\sll_N)$ coincide. 

In this section we show that $\pi^h$ and $\pi^v$ are extended to an action $\ddot{\pi}$ of the quantum toroidal algebra $\tor,$ such that $\pi^h$ is the pull-back of $\ddot{\pi}$ through the homomorphism (\ref{eq:Uh}), and $\pi^v$  is the pull-back of $\ddot{\pi}$ through the homomorphism (\ref{eq:Uv}). The definition of $\ddot{\pi}$ is based on Proposition \ref{p:shift1}.  

Let $\psi_n : \wedge^n\Vaff \rightarrow \wedge^n\Vaff$ be the the map (\ref{eq:psi}) for $\M = (\KK[z^{\pm 1}]\otimes \KK^L)^{\otimes n}.$ That is
\begin{gather} 
\psi_n : z^{m_1}\ef_{a_1}\vf_{\ep_1}\wedge z^{m_2}\ef_{a_2}\vf_{\ep_2}\wedge \cdots \wedge z^{m_n}\ef_{a_n}\vf_{\ep_n} \mapsto \label{eq:psin}\\ 
z^{m_1-\delta_{\ep_1,N}}\ef_{a_1}\vf_{\ep_1+1}\wedge z^{m_2-\delta_{\ep_2,N}}\ef_{a_2}\vf_{\ep_2+1}\wedge \cdots \wedge z^{m_n-\delta_{\ep_n,N}}\ef_{a_n}\vf_{\ep_n+1},  \nonumber
\end{gather} 
where $\vf_{N+1} $ is identified with $\vf_1.$ Let $\F = \oplus_M \F_M.$ We define a semi-infinite analogue $\pinf \in \End(\F)$ of $\psi_n$ as follows. For $m\in \Zint$ we let   
$$ \pinf |-mNL \rangle  = z^{m-1}\ef_1\vf_1\wedge z^{m-1}\ef_2\vf_1\wedge \cdots \wedge z^{m-1}\ef_L\vf_1 \wedge |-mNL \rangle.$$ 
Any vector in $\F$ can be presented in the form $v\wedge |-mNL\rangle ,$ where $v \in \wedge^n\Vaff$ for suitable $n$ and $m.$ Then we set 
$$ \pinf( v\wedge |-mNL\rangle ) = \psi_n(v)\wedge \pinf |-mNL\rangle.$$ 
By using the normal ordering rules it is not difficult to verify that $\pinf$ is well-defined (does not depend on the choice of $m$). Note that $\pinf : \F_M \rightarrow \F_{M+L},$ and that $\pinf$ is invertible. Moreover   
\begin{equation}
\pinfi \pi^h(X_i) \pinf = \pi^h(X_{i-1}) \qquad (i=0,1,\dots,N-1),  \label{eq:cyc}
\end{equation}
where $X=E,F,K$ and the indices are cyclically extended modulo $N.$ 

\begin{propos} \label{p:twist}
For each vector $w \in \F_M$ we have 
\begin{eqnarray}
{\pinfi}\pi^v(\Xt_i(z)){\pinf}(w) & = & \pi^v(\Xt_{i-1}(q^{-1}z))(w), \qquad (2\leq i \leq N-1), \label{eq:tw1} \\
 {\pinfsi}\pi^v(\Xt_1(z)){\pinfs}(w) & = &  \pi^v(\Xt_{N-1}(p^{-1}q^{N-2}z))(w), \label{eq:tw2}
\end{eqnarray}
where $X = E,F,K^{\pm}.$ 
\end{propos}
\begin{proof}
To prove the proposition we use the following lemmas. 

\begin{lemma}\label{l:mlbc}
Let $0 \leq d \leq l,$ $n=s+lNL,$ where $M\equiv s\bmod N,$ $s\in \{0,1,\dots,NL-1\}.$  Let $w= z_1^{\un{k_1}}\ef_{\dot{k}_1}\vf_{\ov{k_1}} \wedge z_2^{\un{k_2}}\ef_{\dot{k}_2}\vf_{\ov{k_2}}  \wedge \cdots \wedge z_n^{\un{k_n}}\ef_{\dot{k}_n}\vf_{\ov{k_n}}$ be a normally ordered wedge from $V_{M,n}^d,$ let $b,c$ be integers such that $1 \leq b \leq L,$ $1 \leq c \leq N$.
We define $f \in \KK[z_1^{\pm 1},\dots,z_{n+bc}^{\pm 1}]\otimes (\KK^L)^{\otimes(n+bc)}$ as follows.

\begin{align}
& f = (z_1^{\un{k_1}}z_2^{\un{k_2}} \cdots z_n^{\un{k_n}})(z_{n+1} \cdots z_{n+bc})^{m}\otimes \label{eq:ffff}\\ & \qquad \qquad \qquad \qquad\otimes (\ef_{\dot{k}_1}\ef_{\dot{k}_2} \cdots \ef_{\dot{k}_n})\underbrace{(\ef_1 \cdots \ef_1)}_{ c \: {\mathrm {times}}}\underbrace{(\ef_2 \cdots \ef_2)}_{ c \:{\mathrm {times}}} \dots \underbrace{(\ef_b \cdots \ef_b)}_{ c \:{\mathrm {times}}}, \nonumber   
\end{align}
where $m = \un{o_{n+1}} =  \un{o_{n+2}} = \cdots =\un{o_{n+bc}}.$ The monomial $f$ given by (\ref{eq:ffff}) satisfies the assumptions of Lemma \ref{l:ml} with $k=bc,$ and $\ov{n}(j) = c$ for all $j \in \{1,2,\dots,b\}.$ Let 
$\KC_{n,bc}^m$ and $\LC_{n,bc}^m$ be the corresponding subspaces of $\KK[z_1^{\pm 1},\dots,z_{n+bc}^{\pm 1}]\otimes (\KK^L)^{\otimes(n+bc)}.$

Let $y= y^{(n)} \otimes ( \vf_{\ep_{1}} \otimes \cdots \otimes \vf_{\ep_{bc}}) \in (\KK^N)^{\otimes n} \otimes (\KK^N)^{\otimes bc} $ such that $N-c+1 \leq \ep_{i} \leq N$ $(1\leq i \leq bc)$, and let $f_1 \in \KC_{n,bc}^m,$ $f_2 \in \LC_{n,bc}^m. $ Then 
\begin{align}
 &\wedge( f_1 \otimes y ) \in \oplus_{d^{\prime} > l} V_{M,n+bc}^{d^{\prime}}, \tag{i} \\ 
 &\wedge( f_2 \otimes y )  = 0. \tag{ii}
\end{align}
\end{lemma}
\begin{proof} \mbox{} \\ 
(i) The vector $\wedge( f_1 \otimes y )$ is a linear combination of normally ordered wedges  
$$ u_{(k_i)} = u_{k_1} \wedge u_{k_2} \wedge \cdots \wedge u_{k_n} \wedge u_{k_{n+1}} \wedge \cdots \wedge u_{k_{n+bc}}  $$ 
such that $\un{k_{n+1}} < \un{o_{n+1}}.$ This inequality implies that $\deg u_{(k_i)} \geq l+1.$ \\
\noindent (ii) 
It is sufficient to show that 
%The normal ordering rules (\ref{eq:n1} -- \ref{eq:n4}) show that the wedge   
%$$ z^m\ef_{a_1}\vf_{\ep_1}\wedge z^m\ef_{a_2}\vf_{\ep_2} \wedge \cdots \wedge z^m\ef_{a_{bc}}\vf_{\ep_{bc}} \in \wedge^{NL}\Vaff $$ 
\begin{equation}
\wedge( \ef_{a_1}\ef_{a_2}\cdots\ef_{a_{bc}} \otimes  \vf_{\ep_{1}} \vf_{\ep_{2}} \cdots \vf_{\ep_{bc}} )\in \wedge^{bc}\Vaff 
\label{lemw}
\end{equation}
is zero whenever there is $J \in \{1,2,\dots,b\}$ such that $\#\{ i \: | \: 1\leq  i\leq bc ,\; a_i = J \} > c.$ 

Using the normal ordering rules (\ref{eq:n1} -- \ref{eq:n4}) one can write (\ref{lemw}) as a linear combination of the normally ordered wedges $\ef_{a^{\prime}_1}\vf_{\ep^{\prime}_1}\wedge \ef_{a^{\prime}_2}\vf_{\ep^{\prime}_2} \wedge \cdots \wedge \ef_{a^{\prime}_{bc}}\vf_{\ep^{\prime}_{bc}}$.

The $\UU_q(\sll_N)$ and $\UU_q(\sll_L)$-weights of the both sides in the normal ordering rules are equal. This implies that $\# \{ i \: | \: a_i ^{\prime } =J \} >c $ and $\# \{ j \: | \:  \exists i , \; \ep^{\prime}_i =j , \;  a_i ^{\prime } =J \} \leq c$. Therefore, there exists some $i$ such that $a_i ^{\prime }= a_{i+1} ^{\prime }$ and $\ep^{\prime}_i= \ep^{\prime}_{i+1}$. 
On the other hand, we know that $\ef_{a_i ^{\prime }}\vf_{\ep^{\prime}_i} \wedge \ef_{a_i ^{\prime }}\vf_{\ep^{\prime}_i} =0$. This implies that  $\wedge( f_2 \otimes y )  = 0.$

%By this assumption, there exists some $k, k^{\prime }$ such that $k \neq k^{\prime }$, $a_k = a_{k^{\prime }}$ and $\ep_k = \ep_{k^{\prime}}$. 
%
%On the other hand by the normal ordering rules (\ref{eq:n1} -- \ref{eq:n4}), we have $\ef_{a_l}\vf_{\ep_l}\wedge \ef_{a_{l^{\prime}}}\vf_{\ep_{l^{\prime}}}= -q^{f(l,l^{\prime })} \ef_{a_{l^{\prime}}}\vf_{\ep_{l^{\prime}}} \wedge\ef_{a_l}\vf_{\ep_l}$, where $f(l,l^{\prime })= -1,0$ or $1$.
%
%By exchanging the wedge, we can set that $\ef_{a_k}\vf_{\ep_k}$ is next to $\ef_{a_{k^{\prime}}}\vf_{\ep_{k^{\prime}}}$, and we know that $\ef_{a_k}\vf_{\ep_k} \wedge \ef_{a_k}\vf_{\ep_k} =0.$
%This implies that  $\wedge( f_2 \otimes y )  = 0.$
\end{proof}

\begin{lemma} \label{l:lt}
Suppose $d$ and $l$ are integers such that $0 \leq d \leq l.$ Let $n= s + lNL,$ where $s \in \{0,1,\dots,NL-1\}$ is defined from $M \equiv s \bmod NL.$ Let $m$ be the integer such that $M-s-lNL = -mNL.$

For $1 \leq b \leq L$ we put 
\begin{align*}
&v_{b,N} = z^m\ef_1\vf_N\wedge z^m\ef_2\vf_N\wedge \cdots \wedge z^m\ef_b\vf_N,\\
&v_{b,N-1} = z^m\ef_1\vf_N\wedge z^m\ef_1\vf_{N-1}\wedge z^m\ef_2\vf_N\wedge z^m\ef_2\vf_{N-1}\wedge \cdots \wedge z^m\ef_b\vf_N\wedge z^m\ef_b\vf_{N-1}.  
\end{align*}

Assume $v \in V^d_{M,s+lNL}.$ Then 
\begin{align}
&\pi^v_{(n+b)}(\Xt_i(z))(v\wedge v_{b,N}) = \pi^v_{(n)}(\Xt_i(z))(v)\wedge  v_{b,N} , \label{eq:lt1}\\
&\pi^v_{(n+2b)}(\Xt_{N-1}(z))(v\wedge v_{b,N-1}) = \pi^v_{(n)}(\Xt_{N-1}(z))(v) \wedge  v_{b,N-1} \label{eq:lt2}
\end{align}
Here $1\leq i \leq N-2.$
\end{lemma}

For the proof, see the appendix.
\vspace{.2in}

%\begin{proof} \mbox{} \\ 
%The proof  of this lemma  closely follows the proof of Lemma 22 in \cite{STU}, to which  the present lemma reduces when  $L=1.$ 
%The only point which is essentially different is that the vanishing condition of $\tilde w$ which appears in the proof of \cite[Lemma 22]{STU} is to be replaced by the $c=2$ or $1$ case of Lemma \ref{l:mlbc} of this article.
%\end{proof}

\noindent Retaining the notations introduced in the statement of the above lemma, we continue the proof of the proposition. We may assume that $w \in \F_M^d.$ Then, by Proposition \ref{p:rho}, $w = v\wedge |-mNL\rangle,$ where $v \in V_{M,s+lNL}^d.$ By Definition  \ref{dp:def}, for $ 2 \leq i \leq N-1$ we have    
\begin{equation} \pi^v(\Xt_{i-1}(q^{-1}z))(v\wedge |-mNL\rangle) = \pi_{(n)}^v(\Xt_{i-1}(q^{-1}z))(v )\wedge |-mNL\rangle.  
\label{eq:inftwist1}\end{equation}
The definition of $\pinf$ yields 
$$ |-mNL\rangle = v_{L,N}\wedge \pinfi|-mNL\rangle,$$ 
where $v_{L,N}$ is defined in the statement of Lemma \ref{l:lt}. 
Applying (\ref{eq:lt1}) in this lemma, we have 
$$\pi^v_{(n+L)}(\Xt_{i-1}(q^{-1}z))(v\wedge v_{L,N}) = \pi^v_{(n)}(\Xt_{i-1}(q^{-1}z))(v)\wedge v_{L,N}. 
$$ 
Taking this, and Proposition \ref{p:fintwist} into account, we find that the right-hand side of  (\ref{eq:inftwist1}) equals 
$$
\psi_{n+L}^{-1}\pi^v_{(n+L)}(\Xt_i(z))\psi_{n+L}(v\wedge v_{L,N})\wedge \pinfi|-mNL\rangle,  
$$
which in turn is equal, by definition of $\pinf,$  to 
\begin{equation}
\pinfi \left(\pi^v_{(n+L)}(\Xt_i(z))\psi_{n+L}(v\wedge v_{L,N})\wedge |-mNL\rangle \right).   
\label{eq:iftwist2}\end{equation}
It is clear, that $\psi_{n+L}(v\wedge v_{L,N}) \in V^{d^{\prime}}_{M+L,n+L}$ for some non-negative integer $d^{\prime}.$ Choosing now $m$ large enough, or, equivalently, $l$ large enough (cf. the statement of Lemma \ref{l:lt}), we have by Definition  \ref{dp:def}:
$$ \pi^v_{(n+L)}(\Xt_i(z))\psi_{n+L}(v\wedge v_{L,N})\wedge |-mNL\rangle  =
 \pi^v(\Xt_i(z))\left(\psi_{n+L}(v\wedge v_{L,N})\wedge |-mNL\rangle \right) .$$
Since $\pinf(v\wedge |-mNL\rangle ) = \psi_{n+L}(v\wedge v_{L,N})\wedge |-mNL\rangle,$ we find that (\ref{eq:iftwist2}) equals 
$$ \pinfi \pi^v(\Xt_i(z))\pinf\left(v\wedge|-mNL\rangle \right).$$
Thus (\ref{eq:tw1}) is proved. 

A proof of (\ref{eq:tw2}) is similar. Here the essential ingredients are the relation (\ref{eq:lt2}), and those relations of Proposition \ref{p:fintwist} which contain the square of $\psi.$   
\end{proof}
Now by Propositions \ref{p:shift1} and \ref{p:twist} we obtain 
\begin{theor}
The following map extends to a representation of $\tor$ on $\F_M.$
\begin{alignat}{4}
&  \ddot{\pi} : X_i(z)  &\quad \mapsto\quad   & \pi^v(\Xt_i(d^iz)) \qquad  (1\leq i <N), \label{eq:tt1}\\   
&  \ddot{\pi} : X_0(z)  &\quad  \mapsto \quad  & {\pinfi}\pi^v(\Xt_1(qd^{-1}z)){\pinf}, \label{eq:tt2} \\
&  \ddot{\pi} : \dc   & \quad \mapsto \quad  & d 1, \\
&  \ddot{\pi} : q^{\frac12 c}  & \quad \mapsto\quad   &  1.  
\end{alignat}
Here $d = p^{-1/N}q ,$ and $X = E,F,K^{\pm}.$
\end{theor}
\mbox{} From (\ref{eq:tt1}) it follows that the vertical (level $0$) $\UNp$-action $\pi^v$ is the pull-back of $\ddot{\pi}$ through the homomorphism (\ref{eq:Uv}). Whereas  from (\ref{eq:tt2}) and (\ref{eq:cyc}) it follows that the horizontal (level $L$) $\UNp$-action $\pi^h$ the pull-back of $\ddot{\pi}$ through the homomorphism (\ref{eq:Uh}). Thus as an $\tor$-module the Fock space $\F_M$ has level $(0,L)$ (cf. Section \ref{sec:tor}).

\subsection{The actions of the quantum toroidal algebra on irreducible integrable highest weight modules of $\UU_q^{\prime}(\agl_N)$ }
Let $\Lambda$ be a level $L$ dominant  integral weight of $\UNp.$ In this section we define an action of the quantum toroidal algebra $\tor$ on the irreducible module   
\begin{equation}
\widetilde{V}(\Lambda) = \KK[H_-]\otimes V(\Lambda) 
\end{equation}
of the algebra $\UU_q^{\prime}(\agl_N) = H\otimes \UNp.$ Here (cf. Section \ref{sec:decomp}) $\KK[H_-]$ is the Fock module of the Heisenberg algebra $H,$ and $V(\Lambda)$ is the irreducible highest weight module of $\UNp$ of highest weight $\Lambda.$  

In Section \ref{sec:torHecke} we defined, for any integral weight $\chi$ of $\sll_L,$ the subalgebra $\UU_q({\mathfrak b}_L)^{\chi}$ of $\ULp.$ A level $N$ action of $\ULp$ on the Fock space $\F_M$ ($M\in \Zint$) was defined in Section \ref{sec:Usl}, so that there is an action $\UU_q({\mathfrak b}_L)^{\chi}$ on $\F_M.$ Recall moreover, that the vertical $\UNp$-action $\pi^v$ on $\F_M,$ and, consequently, the action $\ddot{\pi}$ of $\tor,$ depend on two parameters: $p \in q^{\Zint},$ and $\nu$ which is an integral weight of $\sll_L.$     

\begin{propos}\label{p:inv3}
Suppose $p = q^{-2L},$ and $\nu = - \chi - 2\rho$ for an integral $\sll_L$-weight $\chi.$ Then the action $\ddot{\pi}$ of $\tor$ on $\F_M$ leaves invariant the linear subspace $\UU_q({\mathfrak b}_L)^{\chi}\left(\F_M\right).$   
\end{propos}
\begin{proof}
It is sufficient to prove that both the horizontal $\UNp$-action $\pi^h$ and the vertical $\UNp$-action $\pi^v$ leave $\UU_q({\mathfrak b}_L)^{\chi}\left(\F_M\right)$ invariant. The horizontal action commutes with the action of $\ULp.$ Thus it remains to prove that the vertical action leaves  $\UU_q({\mathfrak b}_L)^{\chi}\left(\F_M\right)$ invariant. Let $w \in \F_M^d$ and let $l \geq d.$ By Proposition \ref{p:rho} there is a unique $v \in V_{M,s+lNL}^d$ such that  
$$ w = v \wedge | M - s - lNL \rangle .$$ 
Here $s\in \{0,1,\dots,NL-1\},$ $M\equiv s\bmod NL.$ \\ 
Let $g$ be one of the generators of $\UU_q({\mathfrak b}_L)^{\chi}$ (cf. \ref{eq:gen}). For all large enough $l$ we have  
\begin{equation}
 g(w) = g(v)\wedge | M - s - lNL \rangle c(g),
\label{eq:inv31}\end{equation} 
where $c(g) = q^{-N}$ if $g = \FL_0,$ and $c(g) = 1$ if $g = \FL_a, \KL_a - q^{\chi(a) - \chi(a+1)} 1$ $(1\leq a < L).$ If $g= \FL_0$ then $g(v) \in V_{M,s+lNL}^{d+1},$ otherwise  $g(v) \in V_{M,s+lNL}^{d}.$    

Let $X$ be an element of $\UNp.$ Provided $l$ is sufficiently large, Definition \ref{dp:def} gives 
$$ \pi^v(X)g(w) = \pi_{(s+lNL)}^v(X)g(v) \wedge | M - s - lNL \rangle c(g) .$$ 
By Proposition \ref{p:inv2} the right-hand side of the last equation is a linear combination of vectors  
\begin{equation} h (v^{\prime})\wedge | M - s - lNL \rangle  , 
\label{eq:inv32}\end{equation}  
where $h$ is again one of the generators of  $\UU_q({\mathfrak b}_L)^{\chi},$ and $v^{\prime}$ belongs to either $V_{M,s+lNL}^{d}$ or $V_{M,s+lNL}^{d+1}.$ Applying (\ref{eq:inv31}) again, the vector  (\ref{eq:inv32}) is seen to be proportional to  
$$ h(v^{\prime}\wedge | M - s - lNL \rangle ).$$ 
Thus the vertical action leaves $\UU_q({\mathfrak b}_L)^{\chi}\left(\F_M\right)$ invariant.
\end{proof}  
Now we use Theorem \ref{t:decofF} to define an action of $\tor$ on $\widetilde{V}(\Lambda).$ Fix the unique $M \in \{0,1,\dots,N-1\}$ such that $\ov{\Lambda} \equiv \ov{\Lambda}_M \bmod \ov{Q}_N.$ Since the dual weights $\dot{\Lambda}^{(M)}$ of $\ULp$ are distinct for distinct $\Lambda,$ from Theorem \ref{t:decofF}
we have the isomorphism of $\UU_q^{\prime}(\agl_N)$-modules:
\begin{equation} \widetilde{V}(\Lambda) \cong \F_M/ \UU_q({\mathfrak b}_L)^{\chi}\left(\F_M\right),\label{eq:is}  
\end{equation}
where $\chi $ is the finite part of $\dot{\Lambda}^{(M)}.$ That is for $\dot{\Lambda}^{(M)} = \sum_{a=0}^{L-1} n_a \dot{\Lambda}_a,$ $\chi =  \sum_{a=1}^{L-1} n_a \dot{\ov{\Lambda}}_a.$ 

By Proposition \ref{p:inv2}, the $\tor$-action $\ddot{\pi}$ with $p= q^{-2L},$ $ \nu = - \chi - 2 \rho,$ factors through the quotient map 
$$ \F_M \rightarrow \F_M/ \UU_q({\mathfrak b}_L)^{\chi}\left(\F_M\right),$$
and therefore by (\ref{eq:is}) induces an action of $\tor$ on $\widetilde{V}(\Lambda).$

\appendix
\section{The proof of lemma \ref{l:lt}}
\setcounter{section}{7}

In this  appendix we  prove Lemma \ref{l:lt}. The idea of the proof is essentially the same as that of the proof of \cite[Lemma 23]{STU}. \\ 
% Lemma 23, but technically we need a little more discussions.
{\bf Lemma \ref{l:lt}.} 
{\em Suppose $d$ and $l$ are integers such that $0 \leq d \leq l.$ Let $n= s + lNL,$ where $s \in \{0,1,\dots,NL-1\}$ is defined from $M \equiv s \bmod NL.$ Let $m$ be the integer such that $M-s-lNL = -mNL.$

For $1 \leq b \leq L$ we put 
\begin{align*}
&v_{b,N} = z^m\ef_1\vf_N\wedge z^m\ef_2\vf_N\wedge \cdots \wedge z^m\ef_b\vf_N,\\
&v_{b,N-1} = z^m\ef_1\vf_N\wedge z^m\ef_1\vf_{N-1}\wedge z^m\ef_2\vf_N\wedge z^m\ef_2\vf_{N-1}\wedge \cdots \wedge z^m\ef_b\vf_N\wedge z^m\ef_b\vf_{N-1}.  
\end{align*}

Assume $v \in V^d_{M,s+lNL}.$ Then 
\begin{align}
&\pi^v_{(n+b)}(\Xt_i(z))(v\wedge v_{b,N}) = \pi^v_{(n)}(\Xt_i(z))(v)\wedge  v_{b,N} , \label{eq:alt1}\\
&\pi^v_{(n+2b)}(\Xt_{N-1}(z))(v\wedge v_{b,N-1}) = \pi^v_{(n)}(\Xt_{N-1}(z))(v) \wedge  v_{b,N-1} \label{eq:alt2}
\end{align}
Here $1\leq i \leq N-2.$}

\begin{proof}
As is mentioned in the proof of Lemma 22 in \cite{STU}, for each $i$ $(1 \leq i \leq N-1),$  the subalgebra of $\UNp$ generated by $\tilde{E}_{i,l'} , \; \tilde{F}_{i,l'}, \; \tilde{H}_{i,m'}, \tilde{K}^{\pm }_{i}$ $(l' \in \Zint , \; m' \in \Zint \setminus \{ 0 \} )$ is in fact generated by only the elements $\tilde{E}_{i,0}, \tilde{F}_{i,0}, \tilde{K}^{\pm }_{i}, \tilde{F}_{i,1} $ and $\tilde{F}_{i,-1}.$ 

By the definition of the representation, every generator of the vertical action $\mbox{U}_v$ preserves the degree in the sense of (\ref{eq:fdeg}).
So it is sufficient to show that the actions of $\tilde{E}_{i,0}, \tilde{F}_{i,0}, \tilde{K}^{\pm }_{i}, \tilde{F}_{i,1} $ and $\tilde{F}_{i,-1}$ satisfy the relations (\ref{eq:alt1}, \ref{eq:alt2}).
 For  $\tilde{E}_{i,0}, \tilde{F}_{i,0}, \tilde{K}^{\pm }_{i}$,
this is shown directly by using the definitions of the actions (\ref{eq:pi1}--\ref{eq:pi3}).
Now we must show that 
\begin{align}
&\pi^v_{(n+b)}(\tilde{F}_{i,\pm 1})(v\wedge v_{b,N}) = \pi^v_{(n)}(\tilde{F}_{i,\pm 1})(v)\wedge v_{b,N}, \label{eq:Fi}\\
&\pi^v_{(n+2b)}(\tilde{F}_{N-1,\pm 1})(v\wedge v_{b,N-1}) = \pi^v_{(n)}(\tilde{F}_{N-1,\pm 1})(v)\wedge v_{b,N-1} \label{eq:FN-1} 
\end{align}
Here $1\leq i \leq N-2.$
\vspace{.15in}

We will prove (\ref{eq:FN-1}).
%In the proof we will use the two different notations: 
%\begin{align}
%& u_{k_1}\wedge u_{k_2} \wedge \cdots \wedge u_{k_{N+2}} \quad \text{ or} \quad  \Lambda ( \zz^{\mm} \otimes {\bold v}_{\ee} ) \\ 
%& \text{ where } \quad k_i = \ep_i - N m_i \quad  \text{ and } \quad \mm = (m_1,\dots,m_{N+2}), \quad \ee  = (\ep_1,\dots,\ep_{N+2})
%\end{align}
%for an element from $\wedge ^{N+2} V(z) = \cplx [ z_{1}^{\pm 1} , \cdots z_{N+2}^{\pm 1}] \otimes (\otimes ^{N+2} V) / \sum_{i=1}^{N+1} $Im$ (\stackrel{c}{T_i} - \stackrel{s}{T_i} )$.  

%The element like $( \; \cdot \; \wedge \; \cdot \; \wedge \; \cdot \; )$ is in $\wedge ^{N+2} V(z)$ and the element like $\Lambda ( \; \cdot \; \otimes \; \cdot \; \otimes \; \cdot \;)$ is in $\cplx [ z_{1}^{\pm 1} , \cdots z_{N+2}^{\pm 1}] \otimes (\otimes ^{N+2} V) / \sum_{i=1}^{N+1} Im (\stackrel{c}{T_i} - \stackrel{s}{T_i} )$.  

%The first ones are $\Lambda (P(\zz ) \otimes w_1 \otimes w_2 )$, where $P(\zz ) \in \cplx [ z_1 ^{\pm 1}, \dots ,z_N' ^{\pm 1}]$, $w_1 \otimes _{i=1}^{N} \in \cplxn$ and $w_2 \in \otimes _{i=N+1}^{N'} \cplxn$.

For any $ M',M'', M''' $ ($1 \leq M', M'', M'''\leq N+2b, \; M'\leq M'' $), we define an $\UNp $--action on the space $ \KK[z_1^{\pm 1},\dots,z_{n+2b}^{\pm 1}]\otimes (\KK^L)^{\otimes(n+2b)} \otimes (\KK^N)^{\otimes (n+2b)}$ in terms of the Chevalley generators as follows:
\begin{align}
& \EN_i( f \otimes \tilde{v} ) =  \sum_{j=M'}^{M''} (q^{-\unun{M'''}}Y_{j}^{(M''')})^{-\delta (i=0)} f \otimes (e_{i,i+1})_{j} (k_{i})_{j+1} \dots  (k_{i})_{M''} \tilde{v}, \label{eM} \\ 
& \FN_i ( f \otimes \tilde{v})    =   \sum_{j=M'}^{M''} (q^{-\unun{M'''}}Y_{j}^{(M''')})^{\delta (i=0)} f \otimes (k_{i}^{-1})_{M'}\dots (k_{i}^{-1})_{j-1}(e_{i+1,i})_{j} \tilde{v}. \label{fM} \\  
& \KN_i ( f \otimes \tilde{v}) =   f \otimes (k_i)_{M'} (k_i)_{M'+1} \dots (k_i)_{M''} \tilde{v}. \label{kM}
%& \ee_0( f \otimes w) =  \sum_{j=M'}^{M''} f \cdot (Y_{j}^{(M''')})^{-1} \otimes E_j^{N,1} K_{j+1}^{0}\dots K_{M''}^{0} w,  \\ 
%& \ff_0 ( f \otimes w)    =   \sum_{j=M'}^{M''} f \cdot Y_{j}^{(M''')}) \otimes (K_{M'}^{0})^{-1}\dots (K_{j-1}^{0})^{-1}E_j^{1,N} w.  \\  
%& \kk _0 ( f\otimes w) = f \otimes (K^0_{M'} K^0_{M'+1} \dots K^0_{M''}) w. \label{kM}   
\end{align}
Here $i= 0, \dots ,N-1$, indices are cyclically extended modulo $N$, $f \in \KK[z_1^{\pm 1},\dots,z_{n+2b}^{\pm 1}]\otimes (\KK^L)^{\otimes(n+2b)}$, $ \tilde{v} \in  (\KK^N)^{\otimes (n+2b)}$, and the meaning of the notations $(e_{i,i'})_j$, $(k ^{\pm 1}_i)_j$ is the same as in Section \ref{sec:VVdual}. It is understood, that for $M''' < n+2b$ the operators  $Y_{i}^{(M''')}$ in (\ref{eM}, \ref{fM}) act non-trivially only on the variables $z_1,z_2,\dots,z_{M'''}$ and on the first $M'''$ factors in $\KK^{\otimes(n+2b)}.$  
Note that the $\UNp$-action is well--defined because of the commutativity of $Y_{i}^{(M''')}$ ($i = 1, \dots , M'''$).
% The other generators act by (\ref{e: Efin}--\ref{e: Kfin}).
 The actions of the Drinfeld generators are determined by the actions of the Chevalley generators.
 
Let $\tilde{X}$ be an element of $\UNp $, we denote by $\tilde{X}^{(M',M''), M'''}$ the operator giving the action of  $\tilde{X}$ on the space $\KK[z_1^{\pm 1},\dots,z_{n+2b}^{\pm 1}]\otimes (\KK^L)^{\otimes(n+2b)} \otimes (\KK^N)^{\otimes (n+2b)}$ in accordance with (\ref{eM}--\ref{kM}).

 Also, we set $\tilde{X}^{\{ j \}, M''' }$ $=$ $\tilde{X}^{(j,j), M'''}$  $(j=1, \dots , M''').$ 

With these definitions, for any two elements $\tilde{X}$ and $\tilde{Y}$ from  $\UNp ,$  the operators  $\tilde{X}^{(M',M''), M'''}$ and $\tilde{Y}^{(N',N''), M'''}$ commute if $M''<N'$ or $N''<M'$. Note that for any $\Xt \in  \UNp$ we have 
$$ \pi^v_{(n+2b)}(\Xt)\wedge(f \otimes \tilde{v}) = \wedge\left( \Xt^{(1,n+2b),n+2b}(f \otimes \tilde{v})\right).$$

%This shows that in this condition if we have $\Delta X= \sum _{\nu } Y_{\nu } \otimes Z_{\nu }$ then $X (f \otimes \tilde{v} ) = \sum _{\nu } Y_{\nu }^{(*)} Z_{\nu }^{(**)} f \otimes \tilde{v} ) $.

\mbox{} 

\noindent Let  $UN_+$ and $UN_-^{2}$ be the left ideals in $\UNp$ generated respectively by $\{ \tilde{E}_{i,k'} \}$ and $\{ \tilde{F}_{i,k'}\tilde{F}_{j,l'} \}$. Let $UN_+^{(M',M''), M'''},(UN_-^2)^{(M',M''), M'''}$ be the images of these ideals with respect to the map $\UNp \rightarrow \End\left(\KK[z_1^{\pm 1},\dots,z_{n+2b}^{\pm 1}]\otimes (\KK^L)^{\otimes(n+2b)} \otimes (\KK^N)^{\otimes (n+2b)} \right)$ given by (\ref{eM}--\ref{kM}). 
Then the following relations hold: 
\begin{align*}
&  \pi_{(n+2b)}^v(\tilde{F}_{N-1,1}) \wedge (f \otimes \tilde{v} ) \equiv \wedge( (\tilde{K}_{N-1} ^{(1,n+2b-2), n+2b}  \tilde{F}_{N-1,1}^{(n+2b-1,n+2b),n+2b}+  \tilde{F}_{N-1,1}^{(1,n+2b-2), n+2b}) (f \otimes \tilde{v} )) , \\
&  \pi_{(n+2b)}^v(\tilde{F}_{N-1,-1}) \wedge (f \otimes \tilde{v} ) \equiv \wedge (((\tilde{K}_{N-1} ^{(1,n+2b-2), n+2b})^{-1}  \tilde{F}_{N-1,-1}^{(n+2b-1,n+2b),n+2b}+  \tilde{F}_{N-1,-1}^{(1,n+2b-2), n+2b}  \\
& \quad \quad \quad + (q^{-1}-q)(\tilde{K}_{N-1}^{(1,n+2b-2), n+2b})^{-1} \tilde{H}_{N-1,-1}^{(1,n+2b-2), n+2b}  \tilde{F}_{N-1,0}^{(n+2b-1,n+2b),n+2b} ) (f \otimes \tilde{v} )) ,\nonumber \\
& \qquad \text{where}\qquad f \in \KK[z_1^{\pm 1},\dots,z_{n+2b}^{\pm 1}]\otimes (\KK^L)^{\otimes(n+2b)}, \quad  \tilde{v} \in (\KK^N)^{\otimes n+2b}. \nonumber
\end{align*}
Here the equivalence $\equiv$ is understood to be modulo $$ \wedge (UN_+^{(1,n+2b-2), n+2b} 
\cdot (UN_-^2)^{(n+2b-1,n+2b),n+2b} (f \otimes \tilde{v} )).$$
These relations follow from the the coproduct formulas which have been 
obtained in \cite[Proposition 3.2.A]{koyama}:
\begin{eqnarray}
& \Delta^+ (\tilde{F}_{i,1}) \equiv \tilde{K}_{i} \otimes \tilde{F}_{i,1} + \tilde{F}_{i,1} \otimes 1 &
\mbox{mod } UN_+ \otimes UN_-^{2} , \label{cop1} \\
& \Delta^+ (\tilde{F}_{i,-1}) \equiv  \tilde{K}_{i}^{-1} \otimes \tilde{F}_{i,-1} + \tilde{F}_{i,-1} \otimes 1 &
\label{cop2} \\
&  + (q^{-1}-q)\tilde{K}_i^{-1} \tilde{H}_{i,-1} \otimes \tilde{F}_{i,0} & \mbox{mod } UN_+
\otimes UN_-^{2}. \nonumber
\end{eqnarray}
Recall the definition of $\Delta^+$ given in (\ref{eq:co1} -- \ref{eq:co4}). 

Let $w= z_1^{\un{k_1}}\ef_{\dot{k}_1}\vf_{\ov{k_1}} \wedge z_2^{\un{k_2}}\ef_{\dot{k}_2}\vf_{\ov{k_2}}  \wedge \cdots \wedge z_n^{\un{k_n}}\ef_{\dot{k}_n}\vf_{\ov{k_n}}$ be a normally ordered wedge from $V_{M,n}^d,$ and define $f \in \KK[z_1^{\pm 1},\dots,z_{n+2b}^{\pm 1}]\otimes (\KK^L)^{\otimes(n+2b)}$ and $\tilde{v} \in (\KK^N)^{\otimes(n+2b)}$ as follows.

\begin{align}
& f = (z_1^{\un{k_1}}z_2^{\un{k_2}} \cdots z_n^{\un{k_n}})(z_{n+1} \cdots z_{n+2b})^{m}\otimes (\ef_{\dot{k}_1}\ef_{\dot{k}_2} \cdots \ef_{\dot{k}_n})(\ef_1 \ef_1 \ef_2 \ef_2 \cdots \ef_b \ef_b), \label{lemf} \\
& \tilde{v} = (\vf_{\ov{k_1}}\vf_{\ov{k_2}} \cdots \vf_{\ov{k_n}})(\vf_N \underbrace{ \vf_{N-1} ) (\vf_N \vf_{N-1}) \dots (\vf_N }_{ b\: {\mathrm {copies}}} \vf_{N-1}),  \label{lemv}
\end{align}
where $m = \un{o_{n+1}} =  \un{o_{n+2}} = \cdots =\un{o_{n+2b}}.$
Then the monomial $f$ satisfies the assumptions of Lemma \ref{l:ml} with $k=2b,$ and $\ov{n}(j) = 2$ for all $j \in \{1,2,\dots,b\}.$
\vspace{.15in}

Now we will show the equality 
\begin{equation}
 \pi_{(n+2b)}^v(\tilde{F}_{N-1,\pm 1})\wedge (f \otimes \tilde{v} )) 
=\wedge(\tilde{F}_{N-1,\pm 1}^{(1,n+2b-2), n+2b} (f \otimes \tilde{v})).
\label{ffstar}
\end{equation}
First let us prove that any element in $UN_+^{(1,n+2b-2), n+2b} \cdot (UN_-^2)^{(n+2b-1,n+2b),n+2b}$ 
annihilates  the vector $f \otimes \tilde{v}$ where $f$ and  $\tilde{v}$ are given by  (\ref{lemf}) and (\ref{lemv}).  
It is enough to show that 
\begin{equation}
 ( \tilde{F}_{i',k'}^{(n+2b-1,n+2b),n+2b}\tilde{F}_{j',l'}^{(n+2b-1,n+2b),n+2b}) ( \bar{v} \otimes (z_{n+2b-1}^m\ef_b \vf_N \otimes z_{n+2b}^m \ef_b \vf_{N-1}))=0,
\end{equation}
for $\bar{v} \in \KK[z_1^{\pm 1},\dots,z_{n+2b-2}^{\pm 1}]\otimes (\KK^L)^{\otimes(n+2b-2)} \otimes (\KK^N)^{\otimes (n+2b-2)}$.
This follows immediately from the observation that $\wt(\vf_N)+\wt( \vf_{N-1})-\ov{\alpha}_{i'}-\ov{\alpha}_{j'}$ is not a $\UU_q(\sll_N)$-weight of $ (\KK^N)^{\otimes 2}$.

Next we will show that $\wedge (\tilde{F}_{N-1,\pm 1}^{(n+2b-1,n+2b),n+2b} (f \otimes \tilde{v}) )= 0$, (here $f$ and $\tilde{v}$ are given by  (\ref{lemf}) and (\ref{lemv})).
By the formulas (\ref{cop1}) and (\ref{cop2}), we have the following identities  modulo 
$\wedge (UN_+^{\{ n+2b-1 \} ,n+2b }  (UN_- ^2) ^{\{ n+2b\},n+2b } ( f \otimes \tilde{v} ))$ (see also \cite{STU}):
\begin{align}
&  \wedge ( \tilde{F}^{(n+2b-1,n+2b),n+2b}_{N-1,1} (f \otimes \tilde{v} )) \equiv \wedge ( (\tilde{K}_{N-1} ^{\{ n+2b-1 \},n+2b }  \tilde{F}_{N-1,1}^{\{  n+2b \},n+2b }+  \tilde{F}_{N-1,1}^{\{ n+2b-1 \},n+2b }) (f\otimes \tilde{v} )) , \\
&  \wedge (\tilde{F}^{(n+2b-1,n+2b),n+2b}_{N-1,-1} (f \otimes \tilde{v} )) \equiv \wedge (((\tilde{K}_{N-1} ^{\{ n+2b-1 \} ,n+2b})^{-1}  \tilde{F}_{N-1,-1}^{\{  n+2b \} ,n+2b}+  \tilde{F}_{N-1,-1}^{\{  n+2b-1 \} ,n+2b}  \\
& \quad \quad \quad + (q^{-1}-q)[ \tilde{E}_{N-1,0}^{\{ n+2b-1 \} ,n+2b}, \tilde{F}_{N-1,-1}^{\{ n+2b-1 \},n+2b } ] \tilde{F}_{N-1,0}^{\{  n+2b \},n+2b } ) (f \otimes \tilde{v} )) ,\nonumber 
\end{align}
%here we used the relation $[ E_{i,0}, F_{i,-1} ]= K_{i}^{-} H_{i,-1} $ which is proved by (\ref{re}). 
The following formula is essentially written in \cite[Proposition  3.2.B]{koyama}:
\begin{align}
& \tilde{F}_{i,\pm 1}^{\{ l\},n+2b } ( f' \otimes (\otimes _{j=1}^{n+2b} \vf_{\ep _j}))  \label{ko2} \\
& = (q^{i-\unun{n+2b}}(Y_l ^{(n+2b)})^{-1})^{\pm 1} f'\otimes (\otimes _{j=1}^{l-1} \vf_{\ep _j}) \otimes 
\delta _{i,\ep _l} \vf_{i+1} \otimes (\otimes _{j=l+1}^{n+2b} \vf_{\ep _j}),
\nonumber
\end{align}
where $f'\in \KK[z_1^{\pm 1},\dots,z_{n+2b}^{\pm 1}]\otimes (\KK^L)^{\otimes(n+2b)} $ and $\otimes _{j=1}^{n+2b} \vf_{\ep _j} \in (\KK^N)^{\otimes (n+2b)}$.

By (\ref{ko2}) we have $(UN_+^{\{ n+2b-1 \},n+2b } (UN_- ^2)^{\{ n+2b \},n+2b } (f \otimes \tilde{v} )) =0 $, and by  (\ref{ko2}) and Lemma \ref{l:ml} we have
\begin{align}
& \tilde{F}_{N-1,\pm 1}^{(n+2b-1,n+2b),n+2b} (f \otimes \tilde{v} )  \equiv \alpha _{\pm 1} \bar{v} \otimes z_{n+2b-1}^m  \ef_b \vf_N \otimes z_{n+2b}^m \ef_b  \vf_N \label{ot}\\
& \qquad \qquad \qquad \qquad
\mbox{ mod } (\KC_{n,2b}^m + \LC_{n,2b}^m ) \otimes (\check{v}^{(n)} \otimes 
(\vf_N \underbrace{ \vf_{N-1} ) \dots (\vf_N }_{ b-1\: {\mathrm {copies}}} \vf_{N-1}) (\vf_N \vf_N)) , \nonumber
\end{align}
Here $ c _{\pm 1}$ are certain coefficients, $ \bar{v} \in \KK[z_1^{\pm 1},\dots,z_{n+2b-2}^{\pm 1}]\otimes (\KK^L)^{\otimes(n+2b-2)} \otimes (\KK^N)^{\otimes (n+2b-2)}$, $\check{v}^{(n)} \in (\KK^N)^{\otimes n}$.

Using the normal ordering rules, we have $\wedge (\bar{v} \otimes z_{n+2b-1}^m  \ef_b \vf_N \otimes z_{n+2b}^m \ef_b \vf_N) =0.$ By  Lemma \ref{l:mlbc}, we have
\begin{equation}
\wedge ((\KC_{n,2b}^m + \LC_{n,2b}^m ) \otimes (\check{v}^{(n)} \otimes 
(\vf_N \underbrace{ \vf_{N-1} ) \dots (\vf_N }_{ b-1\: {\mathrm {copies}}} \vf_{N-1}) (\vf_N \vf_N))) \in \oplus_{d^{\prime} > l} V_{M,n+2b}^{d^{\prime}}.
\label{ott}
\end{equation}
On the other hand the degree of the wedge (\ref{ott}) is equal to deg$(f \otimes \tilde{v} )=d$. Taking into account that $d \leq l,$ we have $\wedge (\tilde{F}_{N-1,\pm 1}^{(n+2b-1,n+2b),n+2b} (f \otimes \tilde{v} ))=0$.

% and $\tilde{w}$ is a linear combination of normally ordered wedges $w(\nn,{\bold \tau}),$ $\nn \in \MC_{s+nl+2}^n,$ ${\bold \tau} \in \EC(\nn)$ (see subsection \ref{sec:lza}) such that for the sequence $\nn = (n_1,n_2, \dots n_{s+nl+2})$ we have    
%\begin{gather}
% n_1,n_2 , \dots ,n_{s+nl+2} \leq  m ,   
%\text{ and } \quad  \#\{n_i | n_i = m \}  <  2. \label{ni2}   
%\end{gather}
%The inequality (\ref{ni2}) implies, in particular, that $n_{s+nl+1} < m.$ 
%By the definition of the space $\wedge ^{N+2} V(z)$, we have $\wedge (c _{\pm 1} z_{n-2b-1}^m z_{n-2b}^m  \bar{v} \otimes \vf_N \otimes \vf_N)=0 $.

%Assume $\tilde{w}\neq 0$. From the  inequality  $n_{s+nl+1} < m$ it follows that $|\tilde{w}|\geq l+1$ (Cf. Proposition 5(ii) in \cite{TU}). On the other hand we have $|\tilde{w}|=|\Lambda (f \otimes w)|=k$. By the condition $k\leq l$ this is a contradiction. Therefore we conclude that $\tilde{w}=0,$ and hence $\Lambda (\tilde{F}_{N-1,\pm 1}^{(n+2b-1,n-2b),n+2b} (f \otimes w))=0$.

Now we prove that $\wedge((\tilde{K}^{(1,n+2b-2),n+2b}_{N-1})^{-1} \tilde{H}^{(1,n+2b-2),n+2b}_{N-1,-1} \tilde{F}_{N-1,0}^{(n+2b-1,n+2b),n+2b} (f \otimes \tilde{v}))$ vanishes. 
We have
\begin{align}
& \wedge ((\tilde{K}^{(1,n+2b-2),n+2b}_{N-1})^{-1} \tilde{H}^{(1,n+2b-2),n+2b}_{N-1,-1} \tilde{F}_{N-1,0}^{(n+2b-1,n+2b),n+2b} \bar{v} \otimes z_{n+2b-1}^m  \ef_b \vf_N \otimes z_{n+2b}^m \ef_b \vf_{N-1}) \\
& = \wedge(( \tilde{K}^{(1,n+2b-2),n+2b}_{N-1})^{-1} \tilde{H}^{(1,n+2b-2),n+2b}_{N-1,-1} \bar{v} \otimes z_{n+2b-1}^m  \ef_b \vf_N \otimes z_{n+2b}^m \ef_b \vf_N) , \nonumber
\end{align}
here $ \bar{v} \in \KK[z_1^{\pm 1},\dots,z_{n+2b-2}^{\pm 1}]\otimes (\KK^L)^{\otimes(n+2b-2)} \otimes (\KK^N)^{\otimes (n+2b-2)}$.
By (\ref{eM} -- \ref{kM}) the operator $\tilde{H}^{(1,n+2b-2),n+2b}_{N-1,-1}$ is a polynomial in the  operators  $(Y_{j}^{(n+2b)})^{\pm 1}$, $(k_l)^{\pm 1}_j,$  $(e_{l,l'})_j $ where $1 \leq j \leq n+2b-2 $ and $ 1 \leq l,l' \leq N$.
By Lemma \ref{l:ml}, we have
%since $z^m v_n \wedge z^m v_{n} = 0$ as implied by the normal ordering rules for $q$-wedges.
\begin{align*}
& (\tilde{K}^{(1,n+2b-2),n+2b}_{N-1})^{-1} \tilde{H}^{(1,n+2b-2),n+2b}_{N-1,-1} \bar{v} \otimes z_{n+2b-1}^m  \ef_b \vf_N \otimes z_{n+2b}^m \ef_b \vf_N \\
& \qquad \qquad \qquad \equiv c(\hat{v} \otimes z_{n+2b-1}^m  \ef_b \vf_N \otimes z_{n+2b}^m \ef_b \vf_N ) \nonumber \\
& \qquad \qquad \qquad \qquad \qquad
\mbox{ mod } (\KC_{n,2b}^m + \LC_{n,2b}^m ) \otimes (\check{v}^{(n)} \otimes 
(\vf_N \underbrace{ \vf_{N-1} ) \dots (\vf_N }_{ b-1\: {\mathrm {copies}}} \vf_{N-1}) (\vf_N \vf_N)) , \nonumber
\end{align*}
Here $ c $ is a certain coefficient, $ \bar{v}, \hat{v} \in \KK[z_1^{\pm 1},\dots,z_{n+2b-2}^{\pm 1}]\otimes (\KK^L)^{\otimes(n+2b-2)} \otimes (\KK^N)^{\otimes (n+2b-2)}$, $\check{v}^{(n)}$ is an element in $(\KK^N)^{\otimes n}$.
 Repeating the arguments given after the relation (\ref{ott}), we have $\wedge((\tilde{K}^{(1,n+2b-2),n+2b}_{N-1})^{-1} \tilde{H}^{(1,n+2b-2),n+2b}_{N-1,-1} \tilde{F}_{N-1,0}^{(n+2b-1,n+2b),n+2b} (f \otimes \tilde{v}))=0.$
Thus we have shown (\ref{ffstar}).
\vspace{.15in}

Repeatedly applying the arguments  that led to (\ref{ffstar}), we have
\begin{equation}
 \pi_{(n+2b)}^v(\tilde{F}_{N-1, \pm 1})(f \otimes \tilde{v})) = \wedge (\tilde{F}^{(1,n),n+2b}_{N-1, \pm 1}(f \otimes \tilde{v})).
\label{nn2b}
\end{equation}

 To prove $\pi_{(n+2b)}^v(\tilde{F}_{N-1, \pm 1})\wedge(f \otimes \tilde{v})) = \wedge (\tilde{F}^{(1,n),n}_{N-1, \pm 1}(f \otimes \tilde{v}))$, we must show that in the right-hand side of (\ref{nn2b}) we can replace $q^{-\unun{n+2b}}Y^{(n+2b)}_i$ by $q^{-\unun{n}}Y^{(n)}_i$ $(1 \leq i \leq n).$

Observe  that $\tilde{F}^{(1,n),n+2b}_{N-1, \pm 1}$ is a polynomial in the operators 
\begin{equation}
(Y_{j}^{(n+2b)})^{\pm 1},\quad (k_l)^{\pm 1}_j,\quad   (e_{l,l'})_j \quad\text{where $1 \leq j \leq n $ and $ 1 \leq l,l' \leq N.$} 
\label{eq:opers}\end{equation}
By Lemma \ref{l:ml} we have
\begin{equation}
(q^{-\unun{n+2b}}Y^{(n+2b)}_i)^{\pm 1} (f \otimes \tilde{v}) \equiv (q^{-\unun{n}}Y^{(n)}_i)^{\pm 1} (f \otimes \tilde{v}) \mbox{ mod } (\KC ^m _{n,2b} + \LC  ^m _{n,2b}) \otimes \tilde{v}.
\end{equation}

For $f^{\prime} \in \KC ^m _{n,2b} + \LC  ^m _{n,2b},$ and $\EC_n$ a polynomial in (\ref{eq:opers}), the vector $f^{\prime}\otimes \EC_n \tilde{v}$
satisfies the assumption of Lemma \ref{l:mlbc}. By this lemma, and by the arguments given after (\ref{ott}), we have 
\begin{equation}
\wedge ( (\KC ^m _{n,2b} + \LC  ^m _{n,2b}) \otimes \EC_n \tilde{v}) =0
\label{lstab}
\end{equation}

Combining (\ref{lstab}), the commutativity of $\EC_n$ and $(Y^{(\tilde{n})}_i)^{\pm 1}$ $(1 \leq i \leq n, \; \tilde{n} = n $ or $n+2b)$, and the fact that $(Y_i^{(\tilde{n})})^{\pm 1} ( \KC ^m _{n,2b} + \LC  ^m _{n,2b}) \subset  ( \KC ^m _{n,2b} + \LC  ^m _{n,2b}),$ we have  $\pi_{(n+2b)}^v(\tilde{F}_{N-1, \pm 1})\wedge (f \otimes \tilde{v})) = \wedge (\tilde{F}^{(1,n),n}_{N-1, \pm 1}(f \otimes \tilde{v}))$. The relation (\ref{eq:FN-1}) follows. 
%
% To prove (\ref{fn-11}) we must show that in the right-hand side of the last 
%equation we can replace $Y^{(n+2b)}_i$ by $Y^{(n)}_i$ $(1 \leq i \leq n).$ By  Lemma \ref{l: lemma3} (\ref{e: ib}), for $1\leq i \leq N$ and a sequence $\mm = (m_1 , m_2, \dots m_{N+2}) \in \Zint ^{N}$ such that $m_1 , m_2 , \dots m_N < m_{N+1} = m_{N+2}=m $ we have:
%\begin{equation}
%\zz^{\mm} (Y_i^{(n+2b)})^{\pm 1}  =  q^{\pm 2} \zz^{\mm} (Y_i^{(n)})^{\pm 1} + [\dots] ,
%\end{equation}
%where $[\dots]$ signifies a linear combination of monomials $\zz^{\nn}$ $\equiv$ $ z_1^{n_1}z_2^{n_2}\dots z_{n+2b}^{n_{n+2b}}$ such that $n_1,n_2,\dots,n_{n+2b}  \leq  m $ and $\#\{n_i | n_i = m \}  <  2$. By the normal ordering rules, we can write
%\begin{align}
%& \wedge (Y^{(n+2b)}_i(z_{n+2b-1}^{m} z_{n+2b}^{m} \bar{v} \otimes (\vf_N \otimes \vf_{N-1}))) \\
%& = \wedge (Y^{(n+2b-2)}_i (z_{n+2b-1}^{m} z_{n+2b}^{m} \bar{v} \otimes (\vf_N \otimes \vf_{N-1}))) + \tilde{w}, \nonumber
%\end{align}
%where  $\bar{v} \in \KK[z_1^{\pm 1},\dots,z_{n+2b-2}^{\pm 1}]\otimes (\KK^L)^{\otimes(n+2b-2)} \otimes (\KK^N)^{\otimes n+2b-2}$ and the $\tilde{w }$ again has the same meaning as the $\tilde{w }$ in  relation (\ref{ot}).  Repeating the discussion after (\ref{ot}) we can show that $\tilde{w}=0$.
%This follows by application of the Lemma \ref{l: lemma3}.
%Hence we get (\ref{eq:FN-1}).
\vspace{.2in}

To prove (\ref{eq:Fi}), consider the tensor product $\KK[z_1^{\pm 1},\dots,z_{n+b}^{\pm 1}]\otimes (\KK^L)^{\otimes(n+b)} \otimes (\KK^N)^{\otimes (n+b)}$, use the formulas (\ref{cop1}), (\ref{cop2}) and continue  the  proof in a way that is  completely analogous to the proof of (\ref{eq:FN-1}). 
%(This case is easier than (\ref{eq:FN-1}).)
\end{proof}
\vspace{.2in}

%------------------------------------------------------------------------------
%References
%------------------------------------------------------------------------------
\newcommand{\BOOK}[6]{\bibitem[{#6}]{#1}{\sc #2}, {\it #3} (#4)#5.}
\newcommand{\JPAPER}[8]{\bibitem[{#8}]{#1}{\sc #2}, `#3',
{\it #4} {\bf #5} (#6) #7.}
\newcommand{\JPAPERS}[9]{\bibitem[{#9}]{#1}{\sc #2}, `#3', {\it #4} #5 #6,
#7 #8.}

\end{document}